\documentclass[12pt]{article}
\usepackage{amsmath}
\usepackage{amssymb}
\usepackage[ctr1,eng]{myarticl}
\usepackage{ifthen}

\newcommand{\Aut }{\mathrm{Aut}}
\newcommand{\bigW }{\widetilde{W}}

\newcommand{\cB }{\mathcal{B}}

\newcommand{\Dchaintwo}[4]{
\rule[-3\unitlength]{0pt}{8\unitlength}
\begin{picture}(14,5)(0,3)
\put(1,2){\ifthenelse{\equal{#1}{l}}{\circle*{2}}{\circle{2}}}
\put(2,2){\line(1,0){10}}
\put(13,2){\ifthenelse{\equal{#1}{r}}{\circle*{2}}{\circle{2}}}
\put(1,5){\makebox[0pt]{\scriptsize #2}}
\put(7,4){\makebox[0pt]{\scriptsize #3}}
\put(13,5){\makebox[0pt]{\scriptsize #4}}
\end{picture}}
\newcommand{\Dchainthree}[6]{
\rule[-3\unitlength]{0pt}{8\unitlength}
\begin{picture}(26,5)(0,3)
\put(1,2){\ifthenelse{\equal{#1}{l}}{\circle*{2}}{\circle{2}}}
\put(2,2){\line(1,0){10}}
\put(13,2){\ifthenelse{\equal{#1}{m}}{\circle*{2}}{\circle{2}}}
\put(14,2){\line(1,0){10}}
\put(25,2){\ifthenelse{\equal{#1}{r}}{\circle*{2}}{\circle{2}}}
\put(1,5){\makebox[0pt]{\scriptsize #2}}
\put(7,4){\makebox[0pt]{\scriptsize #3}}
\put(13,5){\makebox[0pt]{\scriptsize #4}}
\put(19,4){\makebox[0pt]{\scriptsize #5}}
\put(25,5){\makebox[0pt]{\scriptsize #6}}
\end{picture}}
\newcommand{\Dtriangle}[7]{
\rule[-3\unitlength]{0pt}{12\unitlength}
\begin{picture}(18,7)(0,3)
\put(4,4){\ifthenelse{\equal{#1}{l}}{\circle*{2}}{\circle{2}}}
\put(5,4){\line(1,0){8}}
\put(14,4){\ifthenelse{\equal{#1}{r}}{\circle*{2}}{\circle{2}}}
\put(4.4472,4.8944){\line(1,2){4.1056}}
\put(9,14){\ifthenelse{\equal{#1}{t}}{\circle*{2}}{\circle{2}}}
\put(13.5528,4.8944){\line(-1,2){4.1056}}
\put(2,3){\makebox[0pt][r]{\scriptsize #2}}
\put(9,17){\makebox[0pt]{\scriptsize #3}}
\put(16,3){\makebox[0pt][l]{\scriptsize #4}}
\put(6,9){\makebox[0pt][r]{\scriptsize #5}}
\put(12.5,9){\makebox[0pt][l]{\scriptsize #6}}
\put(9,1){\makebox[0pt]{\scriptsize #7}}
\end{picture}}
\newcommand{\cDtriangle}[7]{
\rule[-10\unitlength]{0pt}{\unitlength}
\begin{picture}(18,7)(0,10)
\put(4,4){\ifthenelse{\equal{#1}{l}}{\circle*{2}}{\circle{2}}}
\put(5,4){\line(1,0){8}}
\put(14,4){\ifthenelse{\equal{#1}{r}}{\circle*{2}}{\circle{2}}}
\put(4.4472,4.8944){\line(1,2){4.1056}}
\put(9,14){\ifthenelse{\equal{#1}{t}}{\circle*{2}}{\circle{2}}}
\put(13.5528,4.8944){\line(-1,2){4.1056}}
\put(2,3){\makebox[0pt][r]{\scriptsize #2}}
\put(9,17){\makebox[0pt]{\scriptsize #3}}
\put(16,3){\makebox[0pt][l]{\scriptsize #4}}
\put(6,9){\makebox[0pt][r]{\scriptsize #5}}
\put(12.5,9){\makebox[0pt][l]{\scriptsize #6}}
\put(9,1){\makebox[0pt]{\scriptsize #7}}
\end{picture}}

\newcommand{\End }{\mathrm{End}}
\newcommand{\extWBG}{W^\mathrm{ext}}
\newcommand{\gDd }{\mathcal{D}}
\newcommand{\grgDd }{\mathfrak{D}}

\newcommand{\Hom }{\mathrm{Hom}}

\newcommand{\linc}[2]{#1#2}
\newcommand{\longDchaintwo}[4]{
\rule[-3\unitlength]{0pt}{8\unitlength}
\begin{picture}(18,5)(0,3)
\put(1,2){\ifthenelse{\equal{#1}{l}}{\circle*{2}}{\circle{2}}}
\put(2,2){\line(1,0){14}}
\put(17,2){\ifthenelse{\equal{#1}{r}}{\circle*{2}}{\circle{2}}}
\put(1,5){\makebox[0pt]{\scriptsize #2}}
\put(9,4){\makebox[0pt]{\scriptsize #3}}
\put(17,5){\makebox[0pt]{\scriptsize #4}}
\end{picture}}

\newcommand{\Ndbasis }{\boldsymbol{\mathrm{e}}}

\newcommand{\ndN }{\mathbb{N}}
\newcommand{\ndQ }{\mathbb{Q}}
\newcommand{\ndR }{\mathbb{R}}
\newcommand{\ndZ }{\mathbb{Z}}
\newcommand{\op }{^\mathrm{op}}

\newcommand{\prods}[2]{#1\cdot #2}

\newcommand{\roots }{\boldsymbol{\Delta }}
\newcommand{\scalp }[2]{\langle #1,#2\rangle _{E}}

\newcommand{\WBg }{\mathcal{G}}

\newlength{\mpb}

\title{Classification of arithmetic root systems of rank 3
\thanks{Partially supported by the European Community under a Marie Curie
Intra-European Fellowship}}
\author{I.~Heckenberger}

\begin{document}

\maketitle

\begin{abstract}
The paper purposes to contribute to the classification of pointed Hopf 
algebras by the method of Andruskiewitsch and Schneider.
The structure of arithmetic root systems is enlightened such that their
relation to ordinary root systems associated to semi-simple Lie algebras
becomes more astounding. As an application all arithmetic root
systems of rank 3 are determined. 
The result gives in particular all finite dimensional
rank 3 Nichols algebras of diagonal type over a field of characteristic zero.


\end{abstract}

\section{Introduction}

As demonstrated impressively in the recent paper \cite{a-AndrSchn05},
the method of Andruskiewitsch and Schneider presented in
\cite{a-AndrSchn98} is powerful enough to classify
finite dimensional pointed Hopf algebras such that the coradical is the
group algebra of an abelian group and the base field has characteristic zero.
In the first step of this method one has to determine all finite dimensional
Nichols algebras and describe them in terms of generators and relations.

If the coradical is the group algebra of a finite abelian group and the 
base field is algebraically closed then the Nichols algebra is of diagonal 
type. In \cite{a-Heck04c} to any Nichols algebra of diagonal type
an arithmetic root system and a groupoid are associated which in spirit
are closely related to the root system and the Weyl group of a semi-simple
Lie algebra. This so called Weyl groupoid allows to introduce an
equivalence relation and leaded to the classification of finite dimensional
rank 2 Nichols algebras of diagonal type \cite{a-Heck04e}. Here as a main tool
skew-differential operators on Nichols algebras were used.

The present paper purposes to introduce additional methods in order to
determine all finite dimensional Nichols algebras of diagonal type.
Skew-differential operators could be avoided completely by investigating
intersections of arithmetic root systems with subspaces of $\ndR ^n$.
This allows to find obstructions for the finiteness of an arithmetic
root system using only very elementary calculations.
Additionally, a subgroup of the
Weyl groupoid is introduced which plays a key role for checking the finiteness
of a particular arithmetic root system.
The main result of the paper is Theorem \ref{t-class3} which gives
a complete list of arithmetic root systems of rank 3. The classification
of finite dimensional rank 3 Nichols algebras of diagonal type with help of
this theorem is easily done.

The paper is organized as follows.
In Sections \ref{ssec-found} and \ref{ssec-basres} the main notions are
recalled and introduced and some known results are collected.
In Section \ref{ssec-subsys} it is shown that the intersection of a
subspace of $\ndR ^n$ with an arithmetic root system is again an
arithmetic root system, and a method is given to obtain a basis of this
subsystem.
In Section \ref{ssec-Weyleq} Weyl equivalence is recalled and used to
develop a criterion to check the finiteness of a given arithmetic root
system. In Section \ref{sec-Ars2} the classification result for
rank 2 arithmetic root systems is recalled and some corollaries are stated.
Finally, in Section \ref{sec-Ars3} with Theorem \ref{t-class3} rank 3
arithmetic root systems are classified. This is the largest part of the paper
since one has to consider several subcases.

The classification method in this paper
indicates that there are no principal difficulties anymore to determine
all arithmetic root systems of low rank. To obtain a complete answer in
arbitrary rank will eventually depend on the number of discrete exceptional
examples not related to semi-simple (super) Lie algebras.

In the paper the following notation will be used. If $R$ is a ring, $M$ an
$R$-module, and $S\subset R$, $N\subset M$, then
write $\prods{S}{N}:=\{rm\,|\,r\in S$, $m\in N\}$
and $\linc{S}{N}:=\{\sum _jr_jm_j\,|\,r_j\in S,m_j\in N,r_j=0
\text{ for all but finitely many $j$}\}$. Further,
$\ndN $ denotes the set of natural numbers strictly greater than 0 and
$\ndN _0:=\ndN \cup \{0\}$. The symbol $R_m$, $m\ge 2$, denotes the set of
$m$-th roots of unities in a field $k$.

\section{Arithmetic root systems}
\label{sec-Ars}

\subsection{Foundations}
\label{ssec-found}

In the whole paper let $k$ be a field of characteristic 0
and $k^\ast $ its group of units with respect to multiplication.

Recall that a bicharacter on a group $(G,\cdot )$ (with values in
$k^\ast $) is a map $\chi :G\times G\to k^\ast $
satisfying the equations
\begin{align}
\begin{gathered}
\chi (1,a)=\chi (a,1)=1,\\
\chi (a\cdot b,c)=\chi (a,c)\chi (b,c),\qquad
\chi (a,b\cdot c)=\chi (a,b)\chi (a,c)
\end{gathered}
\end{align}
for all $a,b,c\in G$. Then $\chi \op :G\times G\to k^\ast $,
$\chi \op (a,b):=\chi (b,a)$ is again a bicharacter, and the notation
$\chi \chi \op (a,b):=\chi (a,b)\chi (b,a)$ will be used for $a,b\in G$.

Let $I$ be an index set and $\chi $ a bicharacter on $\ndZ ^I$.
Choose a basis $E=\{\Ndbasis _i\,|\,i\in I\}$ of $\ndZ ^I$.
The numbers $q_{ij}:=\chi (\Ndbasis _i,\Ndbasis _j)$ are called the
\textit{structure constants of} $\chi $ with respect to $E$.
For given $E$ they determine $\chi $ uniquely.

\begin{defin}
Let $E$ be a basis of $\ndZ ^I$, and $\chi $ a bicharacter on $\ndZ ^I$.
Let $q_{ij}$, $i,j\in I$, denote the structure
constants of $\chi $ with respect to $E$.
The \textit{generalized Dynkin diagram}
of the pair $(\chi ,E)$ is a (non-directed) graph $\gDd _{\chi ,E}$
with the following properties:\\
(i) there is a bijective map $\phi $ from $I$ to the vertices of
$\gDd _{\chi ,E}$,\\
(ii) for all $i\in I$ the vertex $\phi (i)$ is labeled by $q_{ii}$,\\
(iii) for all $i,j\in I$ the number $n_{ij}$ of edges between $\phi (i)$
and $\phi (j)$ is either 0 or 1. If $q_{ij}q_{ji}=1$ then $n_{ij}=0$,
otherwise $n_{ij}=1$ and the edge is labeled by $q_{ij}q_{ji}$.
\end{defin}

{}From now on assume that $I=\{1,2,\ldots ,n\}$ for some $n\in \ndN $
and fix a basis $E=\{\Ndbasis _i\,|\,i\in I\}$ of $\ndZ ^n$.
Let also $F$ be a basis of $\ndZ ^n$ and let $\Ndbasis '\in F$.
Suppose that for all $\Ndbasis ''\in F\setminus \{\Ndbasis '\}$
the numbers
\begin{align}\label{eq-m}
\begin{split}
m(\Ndbasis ', \Ndbasis ''):=\min \{m\in \ndN _0\,|\,&
\chi (\Ndbasis ',\Ndbasis ')^m\chi (\Ndbasis ',\Ndbasis '')
\chi (\Ndbasis '',\Ndbasis ')=1\text{ or}\\
&\chi (\Ndbasis ',\Ndbasis ')^{m+1}=1,\ \chi (\Ndbasis ',\Ndbasis ')\not=1\}
\end{split}
\end{align}
exist. Set $m(\Ndbasis ',\Ndbasis ')=-2$. Then the linear map
$s_{\Ndbasis ',F}\in \Aut (\ndZ ^n)$ defined by
\begin{align}
s_{\Ndbasis ',F}(\Ndbasis ''):=\Ndbasis ''+m(\Ndbasis ', \Ndbasis '')
\Ndbasis '
\end{align}
is a reflection, that is $s_{\Ndbasis ',F}^2=\id $ and the rank of
the unique extension of $s_{\Ndbasis ',F}-\id $ to a linear map in
$\End (\ndQ ^n)$ is 1.

Similarly to the notation in \cite{a-Heck04e} let $\bigW $ denote
the groupoid \cite[Sect.~3.3]{b-ClifPres61}
consisting of all pairs $(T,E)$ where $T\in \Aut (\ndZ ^n)$ and
$E$ is a basis of $\ndZ ^n$, and the composition
$(T_1,E_1)\circ (T_2,E_2)$ is defined (and is then equal to
$(T_1T_2,E_2)$) if and only if $T_2(E_2)=E_1$.
Let $E$ be a basis of $\ndZ ^n$ and $\chi $ a bicharacter on $\ndZ ^n$.
Define $W_{\chi ,E}$ as the smallest subgroupoid of $\bigW $
which contains $(\id ,E)$, and if $(\id, F)\in W_{\chi ,E}$ for
a basis $F$ of $\ndZ ^n$ and $\Ndbasis '\in F$, then
$(s_{\Ndbasis ',F},F),(\id ,s_{\Ndbasis ',F}(F))\in W_{\chi ,E}$
whenever $s_{\Ndbasis ',F}$ is defined. It is called the
Weyl groupoid associated to the pair $(\chi ,E)$.
Let $\extWBG _{\chi ,E}$
denote the smallest subgroupoid of $\bigW $ which contains
$W_{\chi ,E}$ and for each basis $F$ of $\ndZ ^n$ with $(\id ,F)\in
W_{\chi ,E}$ and for each $\tau \in \Aut (\ndZ ^n)$ with $\tau (F)=F$
also $(\tau ,F)\in \extWBG _{\chi ,E}$. It is easily checked that
\begin{align}
\extWBG _{\chi ,E}=\{(\tau T,F)\,|\,(T,F)\in W_{\chi ,E},
\tau \in \Aut (\ndZ ^n),\tau T(F)=T(F)\}.
\end{align}

\begin{defin}
The groupoid $W_{\chi ,E}$ respectively $\extWBG _{\chi ,E}$
is called \textit{full}
if $s_{\Ndbasis ',F}$ is well-defined for all bases $F$ of $\ndZ ^n$
with $(\id ,F)\in W_{\chi ,E}$ and for all $\Ndbasis '\in F$.
\end{defin}

Note that $W_{\chi ,E}$ is full if and only if $\extWBG _{\chi ,E}$ is full,
and $W_{\chi ,E}$ is finite if and only if $\extWBG _{\chi ,E}$ is finite.

An arithmetic root system is a triple $(\roots ,\chi ,E)$ such that
$W_{\chi ,E}$ is full and finite and $\roots =\bigcup \{F\,|\,(\id ,F)\in
W_{\chi ,E}\}$. One says that a pair $(\chi ,E)$ respectively an arithmetic
root system $(\roots ,\chi ,E)$ is of \textit{Cartan type}
(see \cite[Definition 5.18]{inp-Andr02})
if for all $\Ndbasis ',\Ndbasis ''\in E$ one has $\chi (\Ndbasis ',\Ndbasis ')
^{m(\Ndbasis ',\Ndbasis '')}\chi \chi \op (\Ndbasis ',\Ndbasis '')=1$.
It is called of \textit{finite Cartan type}
if the matrix $C=(-m(\Ndbasis ',\Ndbasis ''))_{\Ndbasis ',\Ndbasis ''\in E}$
is a Cartan matrix of finite type.

\subsection{Basic results}
\label{ssec-basres}

Let $E=\{\Ndbasis _1,\cdots ,\Ndbasis _n\}$ be a basis of $\ndZ ^n$,
$\chi $ a bicharacter on $\ndZ ^n$,
and $(q_{ij})_{i,j=1,\ldots ,n}$ the set of its structure constants with
respect to $E$. Assume that $(\roots ,\chi ,E)$ is an arithmetic root system.

The definition of $\roots $ and the fact that $|\det s_{\Ndbasis ',F}|=1$
whenever $(\id ,F)\in W_{\chi ,E}$ and $\Ndbasis '\in F$,
give immediately the following.

\begin{lemma}\label{l-multroots}
If $\alpha \in \prods{(\ndN \setminus \{1\})}{\ndZ ^n}$ then $\alpha
\notin \roots $. In particular, if $\alpha \in \roots $ and $r\in \ndR $
then $r\alpha \in \roots $ if and only if $r\in \{1,-1\}$.
\end{lemma}

By \cite[Theorem 1]{a-Heck04e} there is a one-to-one correspondence
between arithmetic root systems and Nichols algebras of diagonal type
with a finite set of restricted Poincar\'e--Birkhoff--Witt generators.
Using this correspondence and \cite[Lemma 14]{a-Rosso98},
\cite[Lemma 3.7(b)]{inp-AndrSchn02} concerning the vanishing of certain
skew-commutators in a Nichols algebra of diagonal type one obtains
the following facts.

\begin{lemma}\label{l-posroots}
(i) For any $i,j\in \{1,\ldots ,n\}$ with $i\not=j$ and for any $m\in \ndN _0$
one has $\Ndbasis _i+m\Ndbasis _j\in \roots $ if and only if $m(\Ndbasis _i,
\Ndbasis _j)\leq m$.\\
(ii) For all bases $F$ of $\ndZ ^n$ such that $(\id ,F)\in W_{\chi ,E}$
one has $\roots =\roots ^+_F\cup -\roots ^+_F$, where $\roots ^+_F=
\roots \cap \linc{\ndN _0}{F}$.
\end{lemma}

In the following the abbrievation $\roots ^+$ will be used for $\roots ^+_E$.

Arithmetic root systems of Cartan type are completely classified.

\begin{thm}\cite[Theorem 4]{a-Heck04c}\label{t-Cartan}
If $(\chi ,E)$ is of Cartan type then $(\roots ,\chi ,E)$ is an arithmetic root
system if and only if $(\chi ,E)$ is of finite Cartan type.
\end{thm}

\subsection{Subsystems of arithmetic root systems}
\label{ssec-subsys}

In this section intersections of hyperplanes with arithmetic root systems
are analyzed.

Fix $n\in \ndN $ and a basis $E=\{\Ndbasis _i\,|\,1\le i\le n\}$ of $\ndZ ^n$.
Let $(\roots ,\chi ,E)$ be an arithmetic root system, $H\subset \ndR ^n$
a hyperplane containing $0$, and set
$\Gamma :=\linc{\ndR }{(\roots \cap H)}\cap \ndZ ^n$
and $\ndR _+:=\{r\in \ndR \,|\,r\ge 0\}$.

\begin{satz}\label{s-subARS}
There exists a unique basis $E_H$ of the subspace
$\linc{\ndR }{(\roots \cap H)}$ of $\ndR ^n$
such that the relations $E_H\subset \roots ^+\cap H\subset
\linc{\ndR _+}{E_H}$ hold.
The triple $(\roots \cap H,\chi |_{\Gamma \times \Gamma },E_H)$ is an
arithmetic root system.
\end{satz}

\begin{bew}
The uniqueness of $E_H$ follows from Lemma \ref{l-multroots}.

Choose $\alpha \in \ndR ^n\setminus \{0\}$ such that $\scalp{\alpha }{H}=0$, where
$\scalp{\cdot }{\cdot }$ is the scalar product on $\ndR ^n$ defined by
$\scalp{\Ndbasis _i}{\Ndbasis _j}=\delta _{ij}$. Set
$\mu:=\min \{|\scalp{\beta }{\alpha }|\,\big|\,\beta \in \roots \setminus H\}$.
Choose a positive number $r\in \ndR $ such that
\begin{align}\label{eq-varepsilon}
\left|\scalp{\varepsilon }{\beta }\right|<\mu /2
\text{ for all $\beta \in \roots $,}
\end{align}
where $\varepsilon :=r\sum _{i=1}^n\Ndbasis _i$.
Then one gets
\begin{align}\label{eq-H}
\begin{split}
|\scalp{\alpha +\varepsilon }{\beta }|=&|\scalp{\varepsilon }{\beta }|>0
\text{ for all } \beta \in \roots \cap H,\\
|\scalp{\alpha +\varepsilon }{\beta }|\ge &|\scalp{\alpha }{\beta }|-
|\scalp{\varepsilon }{\beta }|>\mu /2
\text{ for all }\beta \in \roots \setminus H.
\end{split}
\end{align}
In particular, $\scalp{\alpha +\varepsilon }{\beta }\not=0$ for all
$\beta \in \roots $.
Thus (see the proof of \cite[Theorem 1]{a-Heck04c}) there exists a basis
$F=\{f_1,\ldots ,f_n\}$ of $\ndZ ^n$ such that $F\subset \roots $,
$\scalp{\alpha +\varepsilon }{f}>0$ for all $f\in F$, and $(\roots ,\chi ,F)$
is an arithmetic root system. Without loss of generality
assume that for an $l\in \{0,1,\ldots ,n\}$ one has
$\scalp{\alpha +\varepsilon }{f_i}<\mu /2$ if and only if $i\le l$.
By (\ref{eq-varepsilon}) and (\ref{eq-H})
one has $f_i\in H$ if and only if $i\le l$, and in this case
$f_i\in \roots ^+$ since $\scalp{\varepsilon }{f_i}=
\scalp{\alpha +\varepsilon }{f_i}>0$. Moreover,
if $\beta \in \roots ^+\cap H$ then $0<\scalp{\alpha +\varepsilon }{\beta }<\mu /2$
by (\ref{eq-varepsilon}), and hence $\beta =\sum _{i=1}^lm_if_i$ with
$m_i\in \ndN _0$ for $1\le i\le l$. Set $E_H:=\{f_1,\ldots ,f_l\}$.
\end{bew}

Let $E_H=\{f_1,\ldots ,f_l\}$ be the basis appearing in Proposition
\ref{s-subARS}. In the following the symbols
$\roots (\chi ;f_1,\ldots ,f_l)$ and $\roots ^+(\chi ;f_1,\ldots ,f_l)$
will be used for the set of roots and positive roots, respectively,
of the arithmetic root system
$(\roots \cap H,\chi |_{\Gamma \times \Gamma },E_H)$.

To apply Proposition \ref{s-subARS} to hyperplanes $H$ generated by roots,
the following lemma will be useful.

\begin{lemma}\label{l-lbasis}
Let $(\roots ,\chi ,E)$ be an arithmetic root system of rank $n$ and assume
that $\{f_1,\ldots ,f_l\}\subset \roots $, where $l<n$,
is a linearly independent set over $\ndR $.
Set $H:=\linc{\ndR }{\{f_1,\ldots ,f_l\}}$.
Then $\roots \cap H\subset \linc{\ndR _+}{\{f_1 ,\ldots ,f_l\}}\cup
-\linc{\ndR _+}{\{f_1,\ldots ,f_l\}}$ if and only if
for all $i<l$ and $(m_{i+1},\ldots ,m_l)\in \ndN _0^{l-i}
\setminus \{0^{l-i}\}$ relations
$f_i-\sum _{j=i+1}^lm_jf_j\notin \prods{\ndN}{\roots }$ hold.
\end{lemma}

\begin{bew}
The only if part of the assertion is obvious. The if part will be proven by
induction over $l$, where for $l=1$ the lemma clearly holds.
Using the induction hypothesis and similar
argumentations as in the proof of Proposition \ref{s-subARS}
one can assume that $f_i=\Ndbasis _i$ for all $i\in \{2,\ldots ,l\}$
and $f_1 =m_1\Ndbasis _1+\sum _{j=2}^lm_j\Ndbasis _j$
for certain $m_j\in \ndN _0$, $1\le j\le l$. Since $\{f_1,\ldots ,f_l\}$
is linearly independent, one has also $m_1>0$.
Moreover, $f_1-\sum _{j=2}^lm_jf_j=m_1\Ndbasis _1\in \prods{\ndN }{\roots }$,
and hence $m_j=0$ for all $j\in \{2,\ldots ,l\}$. Thus $f_1=\Ndbasis _1$
by Lemma \ref{l-multroots}, and the claim follows from the
decomposition $\roots =\roots ^+\cup -\roots ^+$ (see Lemma \ref{l-posroots}).
\end{bew}

\subsection{Weyl equivalence}
\label{ssec-Weyleq}

If $\chi $ and $\chi '$ are bicharacters on $\ndZ ^n$ then
the pairs $(\chi ,E)$ and $(\chi ',E)$ are
called \textit{Weyl equivalent} \cite[Definition 2]{a-Heck04e}
if there exists a linear map $T\in \Aut (\ndZ ^n)$, where
$(T,E)\in \extWBG _{\chi ',E}$, such that equation
$\chi (e,e)=\chi '(T(e),T(e))$ holds for all $e\in \ndZ ^n$.
Further, they are called \textit{twist equivalent}
if additionally equation $T(E)=E$ is valid. Two arithmetic root systems
$(\roots ,\chi ,E)$ and $(\roots ',\chi ',E)$ are called
Weyl respectively twist equivalent if the corresponding pairs
$(\chi ,E)$ and $(\chi ',E)$ have this property.

\begin{bsp}\label{ex-Weyleq}
Let $E$ be a basis of $\ndZ ^n$ and $\chi $ a bicharacter on $\ndZ ^n$, where
$n\ge 2$. Choose $i\in \{1,2,\ldots ,n\}$. For all $j\not=i$ set
$m_{ij}:=m(\Ndbasis _i,\Ndbasis _j)$ and
\begin{align}
p_{ij}:=\begin{cases}
1 & \text{if $q_{ii}^mq_{ij}q_{ji}=1$ for some $m\in \ndZ $,}\\
q_{ii}^{-1}q_{ij}q_{ji} & \text{otherwise.}
\end{cases}
\end{align}
By definition, $(\chi ,E)$ is Weyl equivalent to the pair $(\chi ',E)$ such that
$q'_{jl}:=\chi '(\Ndbasis _j,\Ndbasis _l)=
\chi (s_{\Ndbasis _i,E}(\Ndbasis _j),s_{\Ndbasis _i,E}(\Ndbasis _l))$.
As noted in \cite[Equation (3)]{a-Heck04d} the constants $q'_{jl}$ satisfy
the equations
\begin{align}
q'_{ii}=&q_{ii},&q'_{jj}=&p_{ij}^{m_{ij}}q_{jj},&
q'_{ij}q'_{ji}=&p_{ij}^{-2}q_{ij}q_{ji},&
q'_{jl}q'_{lj}=&p_{ij}^{m_{il}}p_{il}^{m_{ij}}q_{jl}q_{lj}
\end{align}
for all $j,l\in \{1,2,\ldots ,n\}\setminus \{i\}$.
\end{bsp}

\begin{bems}
(i) The set of pairs $(\chi ,E)$ can be canonically identified with
braided vector spaces over $k$ of diagonal type,
and with Nichols algebras over $k$ of diagonal type,
respectively. Therefore it makes sense to speak about the
generalized Dynkin diagram of a braided vector space of diagonal type,
a Nichols algebra of diagonal type, and of an arithmetic root system,
respectively.\\
(ii) Two arithmetic root systems have the same
generalized Dynkin diagram if and only if they are twist equivalent.
\end{bems}

\begin{defin}
Let $\chi $ be a bicharacter on $\ndZ ^n$ and $E$ a basis of $\ndZ ^n$.
Let $Q_0$ denote the set of all generalized Dynkin diagrams $\gDd $
such that $\gDd =\gDd _{\chi ',E}$ for a pair $(\chi ',E)$
which is Weyl equivalent to $(\chi ,E)$. For any diagram
$\gDd =\gDd _{\chi ',E}\in Q_0$ and any $\Ndbasis \in E$ such that
$s_{\Ndbasis ,E}$ is defined, put an arrow
$\rho _{\gDd ,\Ndbasis }$ from $\gDd $ to the generalized Dynkin diagram $\gDd '$
of $(\chi ',s_{\Ndbasis ,E}(E))$, and label this arrow by the vertex of $\gDd $
corresponding to $\Ndbasis $.
The oriented graph determined in this way will be denoted by
$\grgDd ^{\chi ,E}$.
\end{defin}

Define
\begin{align}
\begin{split}
\WBg _{\chi ,E}:=\{T\in \Aut (\ndZ ^n)\,|\,&(T,E)\in \extWBG _{\chi ,E},\\
&\chi (T(e),T(e))=\chi (e,e)\text{ for all }e\in \ndZ ^n\}.
\end{split}
\end{align}

\begin{lemma}\label{l-WBg}
The set $\WBg _{\chi ,E}$ is a subgroup of $\Aut (\ndZ ^n)$.
\end{lemma}

\begin{bew}
For a linear map $T\in \Hom (\ndZ ^m,\ndZ ^n)$
write $T_{E'',E'}$ for the matrix of $T$ with respect to the bases
$E'$ of $\ndZ ^m$ and $E''$ of $\ndZ ^n$.

Suppose that $T,T'\in \WBg _{\chi ,E}$ and set $E'=T(E)$. Then by definition
of $\extWBG _{\chi ,E}$, and since $\chi (T(e),T(e))=\chi (e,e)$ for all
$e\in \ndZ ^n$, there exists $T''\in \Aut (\ndZ ^n)$ such that $(T'',E')\in
\extWBG _{\chi ,E}$, $\chi (T''(e),T''(e))=\chi (e,e)$ for all $e\in
\ndZ ^n$, and $T''_{E',E'}=T'_{E,E}$. Then $(T''T,E)\in \extWBG _{\chi ,E}$
and $(T''T)_{E,E}=\id _{E,E'}T''_{E',E'}T_{E',E}=T_{E,E}T'_{E,E}=(TT')_{E,E}$,
and hence $TT'\in \extWBG _{\chi ,E}$.
Similarly one shows that $T^{-1}\in \extWBG _{\chi ,E}$.
\end{bew}

\begin{bem}\label{rem-grgDd}
(i) Let $\chi $ be a bicharacter on $\ndZ ^n$.
If for all $i,j\in \{1,2,\ldots ,n\}$ there exist numbers $n_{ij}\in \ndN $
such that $(q_{ij}q_{ji})^{n_{ij}}=1$ then the graph $\grgDd ^{\chi ,E}$ is
finite. Indeed, then there exists $m\in \ndN $ such that $\chi (e,e)^m=1$
for all $e\in \ndZ ^n$, and hence in Example \ref{ex-Weyleq}
one has $p_{ij}^m=1$ for all $i,j$ with $i\not=j$. Therefore
$\chi '(e,e)^m=1$ for all $e\in \ndZ ^n$, and hence for all pairs $(\chi ',E)$
which are Weyl equivalent to $(\chi ,E)$ equations $\chi '(e,e)^m=1$ are
satisfied for all $e\in \ndZ ^n$.
Now use the fact that there are only finitely many numbers $\zeta $ such that
$\zeta ^m=1$.\\
(ii) If $n=2$ then one can show that
$\grgDd ^{\chi ,E}$ is finite. For bicharacters on $\ndZ ^n$, where $n\ge 3$,
the finiteness of $\grgDd ^{\chi ,E}$ is not known, even under the
restriction that $W_{\chi ,E}$ is full.\\
(iii) If $(\chi ,E)$ and $(\chi ',E)$ are Weyl equivalent then
the groups $\WBg _{\chi ,E}$ and $\WBg _{\chi ',E}$ are isomorphic.
Indeed, similarly to the proof of Lemma \ref{l-WBg} one can show that
if $(T,E)\in \extWBG _{\chi ,E}$ then $S\in \WBg _{\chi ,E}$ if and only if
$S\in \WBg _{\chi ,T(E)}$.
\end{bem}

Choose a well-ordered set $J$ with smallest element 0,
and a subset $W_0=\{(T_j,E)\,|\,j\in J\}$ of
$\extWBG _{\chi ,E}$ containing $(T_0:=\id ,E)$ with the following properties.
\begin{enumerate}
\item
For each $j\in J\setminus \{0\}$ there exist $i\in J$,
$\Ndbasis '\in T_i(E)$, and $\tau \in \Aut (\ndZ ^n)$, such that $i<j$,
$\tau T_i(E)=T_i(E)$, and $T_j=s_{\Ndbasis ',T_i(E)}\tau T_i$.
\item
If $(\chi ',E)$ is Weyl equivalent to $(\chi ,E)$ then
there exists $j\in J$ such that $(\chi \circ (T_j\times T_j),E)$ is twist
equivalent to $(\chi ',E)$.
\item
The pairs $(\chi \circ (T_i\times T_i),E)$ and
$(\chi \circ (T_j\times T_j),E)$ are twist equivalent if and only if $i=j$.
\end{enumerate}
Obviously such a set $W_0$ exists and its cardinality coincides with the
cardinality of the set of vertices of $\grgDd ^{\chi ,E}$.
The assertions of the next two lemmata follow from the definitions of
$\extWBG _{\chi ,E}$,
$\grgDd ^{\chi ,E}$ and $\WBg _{\chi ,E}$ using arguments as in the proof
of Lemma \ref{l-WBg}.

\begin{lemma}\label{l-genWBg}
The group $\WBg _{\chi ,E}$ is generated by the elements
$\tilde{T}_i^{-1}AT_j$, where $i,j\in J$, $E_j=T_j(E)$,
$A\in \{\tau ,s_{\Ndbasis ,E_j}\tau \,|\,\tau \in \Aut (\ndZ ^n),
\tau (E_j)=E_j,\Ndbasis \in E_j\}$,
$\chi (T_i(e),T_i(e))=\chi (AT_j(e),AT_j(e))$ for all $e\in \ndZ ^n$, and
$(\tilde{T}_i^{-1},A(E_j))$
is the unique element in $\extWBG _{\chi ,E}$ satisfying the relation
$(\tilde{T}^{-1}_i)_{AT_j(E),AT_j(E)}=(T^{-1}_i)_{T_i(E),T_i(E)}$.
\end{lemma}

\begin{lemma}\label{l-finextWBG}
Let $\chi $ be a bicharacter on $\ndZ ^n$ and $E$ a basis of $\ndZ ^n$.
Then $W_{\chi ,E}$ is finite if and only if
the graph $\grgDd ^{\chi ,E}$ and the group $\WBg _{\chi ,E}$ are finite.
\end{lemma}

Our long term aim is to determine all finite dimensional Nichols algebras
of diagonal type over an arbitrary base field of characteristic 0.
By \cite[Theorem 1]{a-Heck04e} this is in principle equivalent to the
classification of arithmetic root systems $(\roots ,\chi ,E)$ where $\chi $
takes values in $k^\ast $. The rank 1 case is trivial, and the rank 2 case was
completely solved in \cite{a-Heck04e}. In this paper all
arithmetic root systems of rank 3 are determined. The proof is based
on the classification result in the rank 2 case, Proposition
\ref{s-subARS}, Theorem \ref{t-Cartan}, and Lemmata \ref{l-genWBg} and
\ref{l-finextWBG}.

Recall that a bicharacter $\chi $ on $\ndZ ^I$ (or a braided vector space
$(V,\sigma )$ of diagonal type or a Nichols algebra $\cB (V)$ of diagonal type)
is called \textit{connected} if for all $i,j\in I$ there exists $m\in \ndN _0$
and a sequence $\{i=i_0,i_1,\ldots ,i_{m-1},i_m=j\}$ of elements of $I$
such that the structure constants $q_{ij}$ of $\chi $ (or of $V$, respectively)
satisfy the conditions $q_{i_l,i_{l+1}}q_{i_{l+1},i_l}\not=1$ for all $l\in
\{0,1,\ldots ,m-1\}$. An arithmetic root system $(\roots ,
\chi ,E)$ is called \textit{connected} if $\chi $ is connected.
Note that a pair $(\chi ,E)$
is connected if and only if its generalized Dynkin diagram
$\gDd _{\chi ,E}$ is a connected graph.

Obviously, $I$ can be
decomposed into the disjoint union of its connected components,
$I=\bigcup _\alpha I_\alpha $, and this induces a decomposition
$\roots =\bigcup _\alpha \roots _\alpha $ into disjoint sets $\roots _\alpha $.
This corresponds to the observation in \cite[Lemma 4.2]{a-AndrSchn00}
that a Nichols algebra of diagonal type can be written as the smash product
of its connected components. Therefore only connected arithmetic root systems
respectively Nichols algebras of diagonal type will be considered.
Note that both twist equivalence and Weyl equivalence respect decompositions
into connected components.

\section{Arithmetic root systems of rank 2}
\label{sec-Ars2}

In this section the classification result of rank 2 arithmetic
root systems is recalled and some corollaries are stated.
To do so let $E=\{\Ndbasis _1,\Ndbasis _2\}$ be a basis of $\ndZ ^2$
and $\chi $ a bicharacter on $\ndZ ^2$ (with values in $k^\ast $).
Set $q_{ij}:=\chi (\Ndbasis _i,\Ndbasis _j)$ for all $i,j\in \{1,2\}$.
Then by \cite[Theorem 4]{a-Heck04e} arithmetic root systems
(up to twist equivalence) correspond to generalized Dynkin diagrams appearing
in Table 1. Moreover, two such arithmetic root systems are
Weyl equivalent if and only if their generalized Dynkin diagrams appear
in the same row of Table 1 and can be presented with the same set of fixed parameters.

\setlength{\unitlength}{1mm}
\settowidth{\mpb}{$q_0\in k^\ast \setminus \{-1,1\}$,}

\begin{table}
\begin{tabular}{r|l|l}
 & \text{generalized Dynkin diagrams} & \text{fixed parameters} \\
\hline \hline
 1 &
 \rule[-3\unitlength]{0pt}{8\unitlength}
 \begin{picture}(14,5)(0,3)
 \put(1,2){\circle{2}}
 \put(13,2){\circle{2}}
 \put(1,5){\makebox[0pt]{\scriptsize $q$}}
 \put(13,5){\makebox[0pt]{\scriptsize $r$}}
 \end{picture}
& $q,r\in k^\ast $ \\
\hline
 2 & \Dchaintwo{}{$q$}{$q^{-1}$}{$q$} & $q\in k^\ast \setminus \{1\}$ \\
\hline
 3 & \Dchaintwo{}{$q$}{$q^{-1}$}{$-1$}\ \Dchaintwo{}{$-1$}{$q$}{$-1$}
& $q\in k^\ast \setminus \{-1,1\}$ \\
\hline
 4 & \Dchaintwo{}{$q$}{$q^{-2}$}{$q^2$} & $q\in k^\ast \setminus \{-1,1\}$ \\
\hline
 \rule[-3\unitlength]{0pt}{10\unitlength}
 5 & \Dchaintwo{}{$q$}{$q^{-2}$}{$-1$}\
\Dchaintwo{}{$-q^{-1}$}{$q^2$}{$-1$} &
 \parbox{\mpb}{$q\in k^\ast \setminus \{-1,1\}$,\\ $q\notin R_4$} \\
\hline
 \rule[-3\unitlength]{0pt}{10\unitlength}
 6 & \Dchaintwo{}{$\zeta $}{$q^{-1}$}{$q$}\
 \Dchaintwo{}{$\zeta $}{$\zeta ^{-1}q$}{\ \ \ \ $\zeta q^{-1}$} &
 \parbox{\mpb}{$\zeta \in R_3$,\\
 $q\in k^\ast \setminus \{1,\zeta ,\zeta ^2\}$} \\
\hline
 7 & \Dchaintwo{}{$\zeta $}{$-\zeta $}{$-1$}\
 \Dchaintwo{}{$\zeta ^{-1}$}{$-\zeta ^{-1}$}{$-1$} & $\zeta \in R_3$ \\
\hline
 8 & \Dchaintwo{}{$-\zeta ^{-2}$}{$-\zeta ^3$}{$-\zeta ^2$}\ \
 \Dchaintwo{}{$-\zeta ^{-2}$}{$\zeta ^{-1}$}{$-1$}\
 \Dchaintwo{}{$-\zeta ^2$}{$-\zeta $}{$-1$}
& $\zeta \in R_{12}$ \\
 & \Dchaintwo{}{$-\zeta ^3$}{$\zeta $}{$-1$}\
\Dchaintwo{}{$-\zeta ^3$}{$-\zeta ^{-1}$}{$-1$} & \\
\hline
 9 & \Dchaintwo{}{$-\zeta ^2$}{$\zeta $}{$-\zeta ^2$}\
 \Dchaintwo{}{$-\zeta ^2$}{$\zeta ^3$}{$-1$}\
 \Dchaintwo{}{$-\zeta ^{-1}$}{$-\zeta ^3$}{$-1$}
 & $\zeta \in R_{12}$ \\
\hline
10 & \Dchaintwo{}{$-\zeta $}{$\zeta ^{-2}$}{$\zeta ^3$}\
 \Dchaintwo{}{$\zeta ^3$}{$\zeta ^{-1}$}{$-1$}\
 \Dchaintwo{}{$-\zeta ^2$}{$\zeta $}{$-1$}
 & $\zeta \in R_9$\\
\hline
 \rule[-3\unitlength]{0pt}{10\unitlength}
 11 & \Dchaintwo{}{$q$}{$q^{-3}$}{$q^3$}
 & \parbox{\mpb}{$q\in k^\ast \setminus \{-1,1\}$,\\
 $q\notin R_3$} \\
\hline
 12 & \Dchaintwo{}{$\zeta ^2$}{$\zeta $}{$\zeta ^{-1}$}\
 \Dchaintwo{}{$\zeta ^2$}{$-\zeta ^{-1}$}{$-1$}\
 \Dchaintwo{}{$\zeta $}{$-\zeta $}{$-1$}
 & $\zeta \in R_8$ \\
\hline
 13 & \Dchaintwo{}{$\zeta ^6$}{$-\zeta ^{-1}$}{\ \ $-\zeta ^{-4}$}\ \
 \Dchaintwo{}{$\zeta ^6$}{$\zeta $}{$\zeta ^{-1}$}\ \
 \Dchaintwo{}{$-\zeta ^{-4}$}{$\zeta ^5$}{$-1$}\
 \Dchaintwo{}{$\zeta $}{$\zeta ^{-5}$}{$-1$} &
 $\zeta \in R_{24}$\\
\hline
 14 & \Dchaintwo{}{$\zeta $}{$\zeta ^2$}{$-1$}\
 \Dchaintwo{}{$-\zeta ^{-2}$}{$\zeta ^{-2}$}{$-1$}
 & $\zeta \in R_5$ \\
\hline
 15 & \Dchaintwo{}{$\zeta $}{$\zeta ^{-3}$}{$-1$}\
 \Dchaintwo{}{$-\zeta $}{$-\zeta ^{-3}$}{$-1$}\
 \Dchaintwo{}{$-\zeta ^{-2}$}{$\zeta ^3$}{$-1$}\
 \Dchaintwo{}{$-\zeta ^{-2}$}{$-\zeta ^3$}{$-1$}
 & $\zeta \in R_{20}$\\
\hline
 16 & \Dchaintwo{}{$-\zeta $}{$-\zeta ^{-3}$}{$\zeta ^5$}\
 \Dchaintwo{}{$\zeta ^3$}{$-\zeta ^4$}{$-\zeta ^{-4}$}\
 \Dchaintwo{}{$\zeta ^5$}{$-\zeta ^{-2}$}{$-1$}\
 \Dchaintwo{}{$\zeta ^3$}{$-\zeta ^2$}{$-1$}
 & $\zeta \in R_{15}$\\
\hline
 17 & \Dchaintwo{}{$-\zeta $}{$-\zeta ^{-3}$}{$-1$}\
 \Dchaintwo{}{$-\zeta ^{-2}$}{$-\zeta ^3$}{$-1$}
 & $\zeta \in R_7$\\
\hline
\end{tabular}
\caption{Weyl equivalence for rank 2 arithmetic root systems}
\end{table}

Later some general facts about rank 2 arithmetic root systems
will be needed which can be obtained directly from Table 1.

\begin{satz}\label{s-genrank2}
Let $E=\{\Ndbasis _1,\Ndbasis _2\}$ be a basis of $\ndZ ^2$
and $\chi $ a bicharacter on $\ndZ ^2$.
Set $q_{ij}:=\chi (\Ndbasis _i,\Ndbasis _j)$ for all $i,j\in \{1,2\}$.
If $W_{\chi ,E}$ is full and finite then the following assertions hold.\\
(i) Either equation
$(q_{12}q_{21}-1)(q_{11}+1)(q_{22}+1)(q_{11}q_{12}q_{21}-1)
(q_{12}q_{21}q_{22}-1)=0$ is satisfied
or one has $q_{11}q_{12}^2q_{21}^2q_{22}=-1$, and in the latter case
also one of the relations $q_{11}\in R_3$, $q_{22}\in R_3$ is fulfilled.\\
(ii) If the arithmetic root system $(\roots ,\chi ,E)$ determined by $\chi $ and $E$
is connected then it is Weyl equivalent to an arithmetic root system
$(\roots ',\chi ',E)$ such that
the constants $q'_{ij}:=\chi '(\Ndbasis _i,\Ndbasis _j)$, $i,j\in \{1,2\}$, satisfy
either equation $q'_{22}=-1$ or equations
$q'_{12}q'_{21}q'_{22}=1$, ${q'_{11}}^aq'_{12}q'_{21}=1$, where
$a\in \{1,2,3\}$, or relations $q'_{12}q'_{21}q'_{22}=1$, $q'_{11}\in R_3$.
\end{satz}

Inserting special values for the structure constants $q_{ij}$,
$i,j\in \{1,2\}$, one obtains the following assertions from Table 1.

\begin{lemma}\label{l-specrank2}
Under the conditions of Proposition \ref{s-genrank2} the following
assertions are true.\\
(i) If $q_{12}q_{21}=-1$ and $q_{11}\not=-1$ then $q_{22}=-1$ and
$q_{11}\in R_3\cup R_4\cup R_6$.\\
(ii) If $q_{12}q_{21}=-q_{11}$ and $q_{22}=-1$ then $q_{11}\in R_2\cup R_3
\cup R_4 \cup R_6\cup R_8\cup R_{10}$.\\
(iii) If $q_{12}q_{21}=-q_{11}$ and $q_{22}=-q_{11}^{-1}$ then $q_{11}\in
R_2\cup R_3\cup R_4\cup R_6\cup R_8$.
\end{lemma}

\begin{lemma}\label{l-specrank2-2}
If $q,r\in k^\ast $, $q,r\not=1$, and $(\roots _1,\chi _1,E)$ and
$(\roots _2,\chi _2,E)$ are arithmetic root systems with generalized
Dynkin diagrams
\Dchaintwo{}{$q$}{$r$}{$r^{-1}$} and
\Dchaintwo{}{$q$}{$r$}{$-1$}~, respectively,
then either $q\in R_2\cup R_3$ or $q^mr=1$ with $m\in \{1,2,3\}$.
\end{lemma}

\begin{bew}
If $(\chi _1,E)$ is of Cartan type then the claim of the lemma follows
from Theorem \ref{t-Cartan}. Otherwise $q\in R_{m+1}$ and $r^{m+1}\not=1$
with $m\in \{1,2,3,4\}$. If $m=4$ then $r\in R_{30}$ and $q=-r^{-3}$,
and if $m=3$ then either $r\in R_8,q=r^2$ or $r\in R_{24},q=r^6$.
In all three cases $(\chi _2,E)$ does not give an arithmetic root system.
\end{bew}

\section{Arithmetic root systems of rank 3}
\label{sec-Ars3}

In this section the main result of this paper is formulated and proved.

\begin{thm}\label{t-class3}
Let $k$ be a field of characteristic 0.
Then twist equivalence classes of connected arithmetic root systems
of rank 3 are in one-to-one correspondence to generalized Dynkin diagrams
appearing in Table 2.
Moreover, two such arithmetic root systems are
Weyl equivalent if and only if their generalized Dynkin diagrams appear
in the same row of Table 2 and can be presented with the same set of
fixed parameters.
\end{thm}


\setlength{\unitlength}{1mm}
\begin{table}
\begin{tabular}{r|l|l|l}
row & gener. Dynkin diagrams & fixed parameters & $\WBg _{\chi ,E}$\\
\hline \hline
 1 & \Dchainthree{}{$q$}{$q^{-1}$}{$q$}{$q^{-1}$}{$q$}
& $q\in k^\ast \setminus \{1\}$ & $[3,4]$\\
\hline
 2 & \Dchainthree{}{$q^2$}{$q^{-2}$}{$q^2$}{$q^{-2}$}{$q$}
& $q\in k^\ast \setminus \{-1,1\}$ & $[3,4]$\\
\hline
 3 & \Dchainthree{}{$q$}{$q^{-1}$}{$q$}{$q^{-2}$}{$q^2$}
& $q\in k^\ast \setminus \{-1,1\}$ & $[3,4]$\\
\hline
 4 & \Dchainthree{}{$-1$}{$q^{-1}$}{$q$}{$q^{-1}$}{$q$}\
 \Dchainthree{}{$-1$}{$q$}{$-1$}{$q^{-1}$}{$q$}
& $q\in k^\ast \setminus \{-1,1\}$ & $\ndZ _2\times [3]$ \\
\hline
 5 & \Dchainthree{}{$-1$}{$q^{-2}$}{$q^2$}{$q^{-2}$}{$q$}\
 \Dchainthree{}{$-1$}{$q^2$}{$-1$}{$q^{-2}$}{$q$}
& $q\in k^\ast \setminus \{-1,1\}$ & $\ndZ _2\times [4]$ \\
 & \Dchainthree{}{$q^2$}{$q^{-2}$}{$-1$}{$q^2$}{$-q^{-1}$}
& $q\notin R_4$ &\\
\hline
 6 & \Dchainthree{}{$-1$}{$q^{-1}$}{$q$}{$q^{-2}$}{$q^2$}
& $q\in k^\ast \setminus \{-1,1\}$ & $\ndZ _2\times [4]$ \\
 & \Dchainthree{}{$-1$}{$q$}{$-1$}{$q^{-2}$}{$q^2$}\
 \Dtriangle{}{$q$}{$-1$}{$-1$}{$q^{-1}$}{$q^2$}{$q^{-1}$}
& & \\
\hline
 7 & \Dchainthree{}{$-1$}{$q^{-1}$}{$q$}{$q^{-3}$}{$q^3$}
& $q\in k^\ast \setminus \{-1,1\}$ & $\ndZ _2\times [6]$ \\
 & \Dchainthree{}{$-1$}{$q$}{$-1$}{$q^{-3}$}{$q^3$}\
 \Dtriangle{}{$q$}{$-1$}{$-1$}{$q^{-1}$}{$q^3$}{$q^{-2}$}
& $q\notin R_3$ &\\
 & \Dchainthree{}{$q^3$}{$q^{-3}$}{$-1$}{$q^2$}{$-q^{-1}$} & & \\
\hline
 8 & \Dchainthree{}{$q$}{$q^{-1}$}{$-1$}{$q$}{$q^{-1}$}\
 \Dchainthree{}{$-1$}{$q$}{$-1$}{$q^{-1}$}{$-1$}
& $q\in k^\ast \setminus \{-1,1\}$ & $\ndZ _2^3$ \\
 & \Dchainthree{}{$-1$}{$q^{-1}$}{$q$}{$q^{-1}$}{$-1$}\
 \Dchainthree{}{$-1$}{$q$}{$q^{-1}$}{$q$}{$-1$}
& & \\
\hline 
 9 & \Dchainthree{}{$q$}{$q^{-1}$}{$-1$}{$r^{-1}$}{$r$}
& $q,r,s\in k^\ast \setminus \{1\}$ & $\ndZ _2^3$\\
& \Dchainthree{}{$q$}{$q^{-1}$}{$-1$}{$s^{-1}$}{$s$}
\Dtriangle{}{$-1$}{$-1$}{$-1$}{$q$}{$r$}{$s$} & $qrs=1$, $q\not=r$, & \\
& \Dchainthree{}{$r$}{$r^{-1}$}{$-1$}{$s^{-1}$}{$s$}
& $q\not=s$, $r\not=s$ & \\
\hline
10 & \Dchainthree{}{$q$}{$q^{-1}$}{$-1$}{$q^{-1}$}{$q$}
 & $q\in k^\ast \setminus \{-1,1\}$ & $\ndZ _2\times [4]$\\
& \Dchainthree{}{$q$}{$q^{-1}$}{$-1$}{$q^2$}{$q^{-2}$}
\Dtriangle{}{$-1$}{$-1$}{$-1$}{$q$}{$q$}{$q^{-2}$} & $q\notin R_3$ & \\
\hline
\end{tabular}
\caption{Weyl equivalence for connected rank 3 arithmetic root systems}
\end{table}

\addtocounter{table}{-1}

\begin{table}
\begin{tabular}{r|l|l|l}
row & gener. Dynkin diagrams & fixed param. & $\WBg _{\chi ,E}$\\
\hline \hline
11 & \Dchainthree{}{$\zeta $}{$\zeta ^{-1}$}{$-1$}{$\zeta ^{-1}$}{$\zeta $}
& $\zeta \in R_3$ & $[3,4]$\\
& \rule{24\unitlength}{0pt}\
\Dtriangle{}{$-1$}{$-1$}{$-1$}{$\zeta $}{$\zeta $}{$\zeta $} & & \\
\hline
12 & \Dchainthree{}{$-\zeta ^{-1}$}{$-\zeta $}{$-\zeta ^{-1}$}{$-\zeta $}{$\zeta $} &
$\zeta \in R_3$ & $[3,4]$\\
\hline
13 & \Dchainthree{}{$\zeta $}{$\zeta ^{-1}$}{$\zeta $}{$\zeta ^{-2}$}{$-1$}\
 \Dchainthree{}{$\zeta $}{$\zeta ^{-1}$}{$-\zeta ^{-1}$}{$\zeta ^2$}{$-1$}
& $\zeta \in R_3\cup R_6$ & $[3,4]$ \\
\hline
14 & \Dchainthree{}{$-1$}{$-\zeta $}{$-\zeta ^{-1}$}{$-\zeta $}{$\zeta $}\
 \Dchainthree{}{$-1$}{$-\zeta ^{-1}$}{$-1$}{$-\zeta $}{$\zeta $}
& $\zeta \in R_3$ & $\ndZ _2\times [4]$ \\
 & \Dchainthree{}{$-\zeta ^{-1}$}{$-\zeta $}{$-1$}{$-\zeta ^{-1}$}{$\zeta ^{-1}$} & & \\
\hline
15 & \Dchainthree{}{$-1$}{$\zeta ^{-1}$}{$\zeta $}{$\zeta $}{$-1$} &
$\zeta \in R_3$ & $\ndZ _2\times [6]$ \\
& \Dchainthree{}{$-1$}{$\zeta $}{$-1$}{$\zeta $}{$-1$}\
\Dtriangle{}{$\zeta $}{$-1$}{$\zeta $}{$\zeta ^{-1}$}{$\zeta ^{-1}$}{$\zeta ^{-1}$}
 & & \\
 & \Dchainthree{}{$-1$}{$\zeta ^{-1}$}{$-\zeta ^{-1}$}{$\zeta ^{-1}$}{$-1$} & & \\
\hline
16 & \Dchainthree{}{$-1$}{$\zeta ^{-1}$}{$\zeta $}{$-\zeta ^{-1}$}{$-\zeta $}
 & $\zeta \in R_3$ & $[6]$ \\
 & \Dchainthree{}{$-1$}{$\zeta $}{$-1$}{$-\zeta ^{-1}$}{$-\zeta $}\
\Dtriangle{}{$\zeta $}{$-1$}{$-1$}{$\zeta ^{-1}$}{$-\zeta $}{$-1$} & & \\
 & \Dchainthree{}{$\zeta $}{$-1$}{$-1$}{$-\zeta ^{-1}$}{$-\zeta $}\
\Dchainthree{}{$\zeta $}{$-\zeta ^{-1}$}{$-\zeta $}{$-\zeta ^{-1}$}{$-\zeta $} & & \\
\hline
17 & \Dchainthree{}{$-1$}{$-1$}{$-1$}{$\zeta $}{$-1$}
& $\zeta \in R_3$ & $\ndZ _2^3$ \\
 & \Dchainthree{}{$-1$}{$-1$}{$\zeta $}{$\zeta ^{-1}$}{$-1$}\
\Dtriangle{}{$-\zeta $}{$\zeta $}{$-1$}{$-\zeta ^{-1}$}{$\zeta ^{-1}$}{$-\zeta ^{-1}$}
& & \\
 & \Dchainthree{}{$-1$}{$\zeta $}{$-1$}{$-\zeta $}{$-1$}\
\Dchainthree{}{$-1$}{$\zeta $}{$-\zeta $}{$-\zeta ^{-1}$}{$-1$} & & \\
 & \Dchainthree{}{$-1$}{$\zeta ^{-1}$}{$\zeta ^{-1}$}{$-\zeta ^{-1}$}{$-1$} & & \\
 & \Dchainthree{}{$-1$}{$\zeta ^{-1}$}{$\zeta $}{$-\zeta $}{$-1$}\
\Dtriangle{}{$-1$}{$-1$}{$\zeta $}{$-1$}{$\zeta ^{-1}$}{$-\zeta $} & & \\
 & \Dchainthree{}{$-1$}{$-1$}{$-1$}{$-\zeta ^{-1}$}{$\zeta ^{-1}$} & & \\
\hline
18 & \Dchainthree{}{$\zeta $}{$\zeta ^{-1}$}{$\zeta $}{$\zeta ^{-1}$}{$\zeta ^{-3}$}\
 \Dchainthree{}{$\zeta $}{$\zeta ^{-1}$}{$\zeta ^{-4}$}{$\zeta ^4$}{$\zeta ^{-3}$}
& $\zeta \in R_9$ & $[3,4]$ \\
\hline
\end{tabular}
\caption{Weyl equivalence for connected rank 3 arithmetic root systems}
\end{table}

Comparing the entries in Table 2 and Table 1 one obtains the following.

\begin{folg}
The generalized Dynkin diagram of a rank 2 subsystem of a connected
arithmetic root system of rank $n\ge 3$ appears in one of the
rows 1--7 and 11 in Table 1.
\end{folg}

\begin{bew}
As argued in \cite[Theorem 1]{a-Heck04c}, using Weyl equivalence
one can assume that the basis of the rank 2 subsystem is a subset
of the basis of the rank $n$ arithmetic root system.
Thus it suffices to check all subdiagrams of the generalized Dynkin diagrams
appearing in Table 2.
\end{bew}

The remaining part of this paper is devoted to the proof of Theorem
\ref{t-class3}.

Let $E=\{\Ndbasis _1,\Ndbasis _2,\Ndbasis _3\}$ be a fixed basis of $\ndZ ^3$
and $\chi $ a bicharacter on $\ndZ ^3$ with values in $k^\ast $. Set $q_{ij}:=
\chi (\Ndbasis _i,\Ndbasis _j)$. Suppose that
$(\roots ,\chi ,E)$ is a connected arithmetic root system.


\begin{lemma}\label{l-trianglecond}
If $q_{12}q_{21}\not=1$ and $q_{13}q_{31}\not=1$ then either $q_{23}q_{32}=1$
or $q_{12}q_{21}q_{13}q_{31}q_{23}q_{32}=1$.
\end{lemma}

\begin{bew}
Assume that $q_{23}q_{32}\not=1$.
If $(\chi , E)$ is of Cartan type then Theorem \ref{t-Cartan} gives
a contradiction. Thus without loss of generality $q_{11}\in R_{m+1}$
with $m\ge 1$ and $(q_{12}q_{21})^{m+1}\not=1$.
By Lemma \ref{l-posroots}(i) one has $a\Ndbasis _1+\Ndbasis _2\in \roots ^+$
for $0\le a\le m$. By the same reason $\Ndbasis _1+\Ndbasis _3,
\Ndbasis _2+\Ndbasis _3\in \roots ^+$.

By Proposition \ref{s-subARS} $\roots (\chi ;\Ndbasis _1,\Ndbasis _2+\Ndbasis _3)$
is finite, and hence Lemmata \ref{l-lbasis} and \ref{l-posroots}(i) imply that
$(m+1)\Ndbasis _1+\Ndbasis _2+\Ndbasis _3\notin \roots ^+$. Since
$(m+1)\Ndbasis _1+\Ndbasis _2+\Ndbasis _3=
(\Ndbasis _1+\Ndbasis _3)+(m\Ndbasis _1+\Ndbasis _2)$, Lemma \ref{l-posroots}(i)
implies that
\begin{align}\label{eq-trianglecond-1}
1=\chi \chi \op (\Ndbasis _1+\Ndbasis _3,m\Ndbasis _1+\Ndbasis _2)
=&q_{11}^{-2}q_{12}q_{21}(q_{13}q_{31})^mq_{23}q_{32}.
\end{align}
If either $m=1$ or $m\ge 2$, $q_{11}q_{13}q_{31}=1$ then the claim holds.
Otherwise $2\Ndbasis _1+\Ndbasis _3\in \roots ^+$, and
$(m+2)\Ndbasis _1+\Ndbasis _2+\Ndbasis _3\notin \roots ^+$ implies that
\begin{align*}
1=\chi \chi \op (2\Ndbasis _1+\Ndbasis _3,m\Ndbasis _1+\Ndbasis _2)
=q_{11}^{-4}(q_{12}q_{21})^2(q_{13}q_{31})^mq_{23}q_{32}.
\end{align*}
Together with (\ref{eq-trianglecond-1}) one obtains that
$q_{11}^2(q_{12}q_{21})^{-1}=1$ which is a contradiction to the assumption
$(q_{12}q_{21})^{m+1}\not=1$.
\end{bew}

\begin{folg}\label{f-rho}
If $q_{12}q_{21}\not=1$ and $q_{13}q_{31}\not=1$ then $\Ndbasis _1+\Ndbasis _2
+\Ndbasis _3\in \roots ^+$.
\end{folg}

\begin{bew}
By Proposition \ref{s-subARS} the set
$\roots (\chi ;\Ndbasis _3,\Ndbasis _1+\Ndbasis _2)$ is finite. By Lemma
\ref{l-trianglecond} one has $\chi \chi \op (\Ndbasis _3,\Ndbasis _1+\Ndbasis _2)
\not=1$ and hence the claim follows from Lemma \ref{l-posroots}(i).
\end{bew}

\begin{lemma}\label{l-t2t3t7t14}
Suppose that $q_{12}q_{21}\not=1$, $q_{13}q_{31}\not=1$, and
$(2a+1)\Ndbasis _1+2\Ndbasis _2\in \roots ^+$ for some $a\in \ndN _0$. Then
$(2b+1)\Ndbasis _1+2\Ndbasis _3\notin \roots ^+$ for all $b\in \ndN _0$.
In particular, either the generalized Dynkin diagram of
$(\chi |_{\linc{\ndZ }{\{\Ndbasis _1,\Ndbasis _3\}}\times
\linc{\ndZ }{\{\Ndbasis _1,\Ndbasis _3\}}},\{\Ndbasis _1,\Ndbasis _3\})$
appears in one of rows 2--7 of Table 1, or one of the following is true:
\begin{itemize}
\item[(i)] $q_{13}q_{31}\in R_{12}$, $q_{11}=-(q_{13}q_{31})^3$, $q_{33}=-1$,
\item[(ii)] $q_{11}\in R_{12}$, $q_{11}^3q_{13}q_{31}=1$, $q_{33}=-1$,
\item[(iii)] $q_{11}\in R_{18}$, $q_{11}^4q_{13}q_{31}=1$, $q_{33}=-1$.
\end{itemize}
\end{lemma}

\begin{bew}
Set $\alpha _a:=(2a+1)\Ndbasis _1+2\Ndbasis _2$ and assume that
$\beta _b:=(2b+1)\Ndbasis _1+2\Ndbasis _3\in \roots ^+$ for some
$b\in \ndN _0$. Then by Lemma \ref{l-lbasis} and Proposition \ref{s-subARS}
one obtains that
$\{\alpha _a,\beta _b\}$ is a basis of $\roots (\chi ;\alpha _a,\beta _b)$
and hence
\begin{align}\label{eq-t2t3t7t14-1}
(a+b+1)\Ndbasis _1+\Ndbasis _2+\Ndbasis _3=&\alpha _a/2+\beta _b/2
\notin \roots (\chi ;\alpha _a,\beta _b).
\end{align}
Moreover $2(a+b+1)\Ndbasis _1+2\Ndbasis _2+2\Ndbasis _3=\alpha _a+\beta _b
\notin \roots $ by Lemma \ref{l-multroots}.
However $(a+1)\Ndbasis _1+\Ndbasis _2,
(b+1)\Ndbasis _1+\Ndbasis _3\in \roots ^+$, and hence
(again by Proposition \ref{s-subARS} and Lemmata \ref{l-lbasis} and
\ref{l-posroots}(i))
\begin{align*}
\chi \chi \op (a\Ndbasis _1+\Ndbasis _2,(b+1)\Ndbasis _1+\Ndbasis _3)=
\chi \chi \op ((a+1)\Ndbasis _1&+\Ndbasis _2,b\Ndbasis _1+\Ndbasis _3)\\
&=\chi \chi \op (\alpha _a,\beta _b)=1.
\end{align*}
The first equation gives $q_{11}^{2(a-b)}q_{12}q_{21}(q_{13}q_{31})^{-1}=1$,
and the third expression divided by the 4th power of the first one
yields
$q_{11}^{-4a+4b+2}(q_{12}q_{21})^{-2}(q_{13}q_{31})^2=1$.
Thus $q_{11}^2=1$, and hence $a=b=0$. This is a contradiction to
Corollary \ref{f-rho} and relation (\ref{eq-t2t3t7t14-1}).

The last assertion follows from \cite[Theorem 4]{a-Heck04e} and the
Appendix of \cite{a-Heck04a}, from which the structures of the root systems
of all Nichols algebras from Table 1 can be obtained.
\end{bew}

\begin{lemma}\label{l-cycle}
Assume that $q_{12}q_{21}\not=1$, $q_{13}q_{31}\not=1$, and
$q_{23}q_{32}\not=1$. Then $q_{11}=-1$ or $q_{22}=-1$ or $q_{33}=-1$.
\end{lemma}

\begin{bew}
Suppose that the claim is false.
If $(\chi ,E)$ is of Cartan type then $\roots $ is not finite by Theorem
\ref{t-Cartan}. Thus without loss of generality one can assume that
$q_{11}\in R_{m+1}$ and $(q_{12}q_{21})^{m+1}\not=1$ for some $m\ge 2$.
The proof of Lemma \ref{l-trianglecond} shows that $q_{11}q_{13}q_{31}=1$.
Further, relations $q_{33}\not=-1$, $q_{11}\not=-1$ and
$q_{11}q_{12}q_{21}\not=1$
and Lemma \ref{l-t2t3t7t14} imply that $q_{33}q_{23}q_{32}=1$. By Lemma
\ref{l-trianglecond} one obtains that $q_{12}q_{21}=
(q_{13}q_{31}q_{23}q_{32})^{-1}=q_{11}q_{33}$. Again by Lemma \ref{l-t2t3t7t14}
and since $q_{22},q_{33}\not=-1$ one has
$(q_{12}q_{21}q_{22}-1)(q_{13}q_{31}q_{33}-1)=0$. If $q_{13}q_{31}q_{33}=1$
then $q_{33}=(q_{13}q_{31})^{-1}=q_{11}$ and hence $q_{12}q_{21}=q_{11}^2$
which contradicts to the assumptions $q_{11}\in R_{m+1}$,
$(q_{12}q_{21})^{m+1}\not=1$. Therefore $q_{22}=(q_{12}q_{21})^{-1}
=q_{11}^{-1}q_{33}^{-1}$.
By Proposition \ref{s-subARS} the set $\roots (\chi ;\Ndbasis _1+\Ndbasis _2,
\Ndbasis _3)$
\longDchaintwo{}{$q_{11}$}{$q_{11}^{-1}q_{33}^{-1}$}{$q_{33}$}
is finite.
This is a contradiction to Proposition \ref{s-genrank2}(i),
Lemma \ref{l-specrank2}(i), and the conditions
$q_{ii}\not=-1$ for all $i\in \{1,2,3\}$.
\end{bew}

\begin{lemma}\label{l-2forkcond}
If $q_{12}q_{21}\not=1$, $q_{13}q_{31}\not=1$, and $q_{23}q_{32}=1$
then the equation
\begin{align*}
(q_{11}+1)(q_{11}q_{12}q_{21}-1)(q_{11}q_{13}q_{31}-1)
(q_{11}q_{12}q_{21}q_{13}q_{31}+1)=0
\end{align*}
holds.
\end{lemma}

\begin{bew}
Suppose that the claim is false.
By Lemma \ref{l-t2t3t7t14} one can assume that $(2b+1)\Ndbasis _1+2\Ndbasis _3
\notin \roots $ for all $b\in \ndN _0$. In particular, one has
\begin{align}\label{eq-2forkcond-1}
(q_{13}q_{31}q_{33}-1)(q_{33}+1)=0.
\end{align}
Further, if $2\Ndbasis _1+\Ndbasis _2\notin \roots ^+$ or
$2\Ndbasis _1+\Ndbasis _3\notin \roots ^+$ then Lemma
\ref{l-posroots}(i) gives a contradiction.

If $2\Ndbasis _1+\Ndbasis _2+2\Ndbasis _3\notin \roots $ then
$2\Ndbasis _1+\Ndbasis _2+2\Ndbasis _3\notin \roots (\chi ;
\Ndbasis _1+\Ndbasis _3,\Ndbasis _2)$, and hence
$(q_{11}q_{13}q_{31}q_{33}+1)(q_{11}q_{13}q_{31}q_{33}q_{12}q_{21}-1)=0$.
Inserting (\ref{eq-2forkcond-1}) one gets again a contradiction. Therefore
$2\Ndbasis _1+\Ndbasis _2+2\Ndbasis _3\in \roots (\chi ;2\Ndbasis _1+\Ndbasis _2,
\Ndbasis _3)$, and hence $q_{33}\not=-1$. By (\ref{eq-2forkcond-1}) one obtains that
\begin{align}\label{eq-2forkcond-2}
q_{13}q_{31}q_{33}=1,\quad q_{33}^2\not=1.
\end{align}
Under this restriction the second sentence of the proof is by Table 1 equivalent
to the equation
\begin{align}\label{eq-2forkcond-3}
(q_{11}^2-q_{33})(q_{11}^2+q_{11}+1)=0.
\end{align}
Since $\roots (\chi ;\Ndbasis _1,\Ndbasis _2)$ is finite, Proposition
\ref{s-genrank2}(i) implies that
\begin{align}\label{eq-2forkcond-4}
(q_{22}+1)(q_{12}q_{21}q_{22}-1)(q_{11}q_{12}^2q_{21}^2q_{22}+1)=0.
\end{align}

Assume that $2\Ndbasis _1+2\Ndbasis _2+\Ndbasis _3\notin \roots $. Then it is
not in $\roots (\chi ;\Ndbasis _1+\Ndbasis _2,\Ndbasis _3)$, and hence
$(q_{11}q_{12}q_{21}q_{22}+1)(q_{11}q_{12}q_{21}q_{22}q_{13}q_{31}-1)=0$.
If $q_{22}=-1$ or $q_{12}q_{21}q_{22}=1$ then the claim of the Lemma
would hold. If $q_{11}q_{12}^2q_{21}^2q_{22}=-1$ then one gets
$q_{12}q_{21}=-q_{13}q_{31}$. On the other hand,
$2\Ndbasis _1+2\Ndbasis _2+\Ndbasis _3\notin
\roots (\chi ;2\Ndbasis _1+\Ndbasis _3,\Ndbasis _2)$, and hence
Lemma \ref{l-posroots}(i) implies that $(q_{12}^2q_{21}^2-1)
(q_{22}+1)(q_{12}^2q_{21}^2q_{22}-1)=0$. Therefore one obtains again a
contradiction. Thus one has $2\Ndbasis _1+2\Ndbasis _2+\Ndbasis _3
\in \roots ^+(\chi ;2\Ndbasis _1+\Ndbasis _3,\Ndbasis _2)$, and hence
$q_{22}\not=-1$.

According to (\ref{eq-2forkcond-3}) consider first the case $q_{11}\in R_3$.
Then one has $3\Ndbasis _1+(\Ndbasis _1+\Ndbasis _2+\Ndbasis _3)
\notin \roots $. This implies that $\chi \chi \op (2\Ndbasis _1+\Ndbasis _2,
2\Ndbasis _1+\Ndbasis _3)=1$, that is $q_{11}q_{12}q_{21}q_{13}q_{31}=1$.
By (\ref{eq-2forkcond-4}) and since $q_{22}\not=-1$, one has
either $q_{12}q_{21}q_{22}=1$ or $q_{11}q_{12}^2q_{21}^2q_{22}=-1$. In the
first case $\roots (\chi ;\Ndbasis _1+\Ndbasis _2+\Ndbasis _3,\Ndbasis _1)$
is of infinite Cartan type, and in the second it is not finite by
Proposition \ref{s-genrank2}(i).

Finally suppose that $q_{11}^3\not=1$. Then $q_{11}^2=q_{33}$ by
(\ref{eq-2forkcond-3}). Recall that (\ref{eq-2forkcond-4}) holds and
$q_{22}\not=-1$. If $q_{12}q_{21}q_{22}=1$ then $\roots (\chi ;
\Ndbasis _1+\Ndbasis _2+\Ndbasis _3,\Ndbasis _1)$
\longDchaintwo{}{$q_{11}$}{$q_{12}^{-1}q_{21}^{-1}$}{$q_{11}$}
is infinite both
if $q_{11}^mq_{12}q_{21}=1$ for some $m\ge 2$ (Cartan type) and if
$q_{11}\in R_{m+1}$, $m\ge 3$ (Proposition \ref{s-genrank2}(i)).
On the other hand, if $q_{11}q_{12}^2q_{21}^2q_{22}=-1$ then
$\roots (\chi ;\Ndbasis _3,2\Ndbasis _1+\Ndbasis _2)$
\Dchaintwo{}{$q_{11}^2$}{$q_{11}^{-4}$}{$-q_{11}^3$}~
is infinite
by Proposition \ref{s-genrank2}(i).
\end{bew}

\begin{lemma}\label{l-bigleg}
If $(\roots ,\chi ,E)$ is a connected arithmetic root system such that
$q_{12}q_{21}\not=1$, $q_{13}q_{31}\not=1$, $q_{23}q_{32}=1$,
$(q_{13}q_{31}q_{33}-1)(q_{33}+1)\not=0$, and
$(q_{11}+1)(q_{11}q_{13}q_{31}-1)=0$, then
$(q_{11}q_{12}q_{21}q_{13}^2q_{31}^2q_{33}-1)(q_{11}^2q_{12}q_{21}q_{13}q_{31}-1)=0$.
\end{lemma}

\begin{bew}
Consider first the case $q_{11}=-1$.
Since $\Ndbasis _1+2\Ndbasis _3\in \roots ^+$,
one obtains the relation $\Ndbasis _1+\Ndbasis _2+2\Ndbasis _3\in 
\roots (\chi ; \Ndbasis _1+2\Ndbasis _3,\Ndbasis _2)\subset \roots $.
Thus $q_{11}=-1$ implies that
$3\Ndbasis _1+\Ndbasis _2+2\Ndbasis _3\notin \roots (\chi ;
\Ndbasis _1+\Ndbasis _2+2\Ndbasis _3,\Ndbasis _1)$.
Considering $\roots (\chi ;
\Ndbasis _1+\Ndbasis _3,\Ndbasis _1+\Ndbasis _2)$ the latter gives
$(q_{12}q_{21}q_{13}q_{31}-1)(-q_{13}q_{31}q_{33}+1)
(-q_{13}q_{31}q_{33}q_{12}q_{21}q_{13}q_{31}-1)=0$.

Assume now that $q_{11}q_{13}q_{31}=1$. If $q_{22}=-1$
then because of Proposition \ref{s-genrank2}(i) the finiteness
of $\roots (\chi ;\Ndbasis _1+\Ndbasis _2,\Ndbasis _3)$ implies that
$$(-q_{11}q_{12}q_{21}+1)(-q_{11}q_{12}q_{21}q_{13}q_{31}-1)
(-q_{11}q_{12}q_{21}q_{13}^2q_{31}^2q_{33}+1)=0.$$
Therefore by Lemma \ref{l-t2t3t7t14} it remains
to consider the case $q_{12}q_{21}q_{22}=1$, and
one can assume additionally that $(q_{11}q_{12}q_{21}-1)
(q_{11}^2q_{12}q_{21}-1)(q_{11}^2+q_{11}+1)=0$.

If $q_{11}q_{12}q_{21}=1$ then the claim of the lemma is obviously true.
If $q_{11}\in R_3$ then
$3(\Ndbasis _1+\Ndbasis _2)+(\Ndbasis _1+2\Ndbasis _3)\notin \roots ^+$,
and hence $4\Ndbasis _1+3\Ndbasis _2+2\Ndbasis _3\notin
\roots ^+(\chi;\Ndbasis _1+\Ndbasis _2+\Ndbasis _3,2\Ndbasis _1+\Ndbasis _2)$
\longDchaintwo{}{$q_{33}$}{$q_{11}^2q_{22}^{-1}$}{~~~$q_{11}q_{22}^{-1}$}.
Therefore Lemma \ref{l-posroots}(i) gives the relation
$(q_{11}^2q_{12}q_{21}-1)(q_{33}+1)(q_{11}^2q_{12}q_{21}q_{33}-1)=0$.

Finally assume that $q_{11}q_{13}q_{31}=q_{11}^2q_{12}q_{21}=
q_{12}q_{21}q_{22}=1$. Then $(\chi ,E)$ is not of Cartan type by Theorem
\ref{t-Cartan}. Thus by the finiteness of
$\roots (\chi ;\Ndbasis _1,\Ndbasis _3)$ and by Table 1 one has
$q_{33}\in R_{a+1}$ with $a\in \{2,3,4\}$. Note that $\Ndbasis _1+\Ndbasis _2
+2\Ndbasis _3\in \roots (\chi ;\Ndbasis _1+2\Ndbasis _3,\Ndbasis _2)$.
Moreover, relations
$(2\Ndbasis _1+\Ndbasis _2)+b\Ndbasis _3\notin \roots $ hold for all $b>a$,
which in turn imply that $\chi \chi \op (\Ndbasis _1+\Ndbasis _2+m\Ndbasis _3,
\Ndbasis _1+a\Ndbasis _3)=1$ whenever $0<m<a$. Insert $m=1$ and $m=2$.
Since $(q_{13}q_{31})^{a+1}\not=1$,
one gets $a=2$ and $q_{11}^3=q_{33}$, and hence
$q_{11}q_{12}q_{21}q_{13}^2q_{31}^2q_{33}=q_{11}^{-3}q_{33}=1$.
\end{bew}

\begin{lemma}\label{l-specbigleg}
If $(\roots ,\chi ,E)$ is a connected arithmetic root system with generalized
Dynkin diagram
\begin{align*}
\Dchainthree{}{$qr^{-1}$}{$q^{-1}$}{$q$}{$r$}{$s$}
\end{align*}
such that $q,r\not=1$, $q^2\not=r^2$, then $q=-1$ or $qr=1$ or $qr^{-1}\in R_3$.
\end{lemma}

\begin{bew}
One can assume that $q_{11}=q$, $q_{22}=s$, $q_{33}=qr^{-1}$.
Thus by Lemma \ref{l-t2t3t7t14} one has $(rs-1)(s+1)=0$. Further,
Proposition \ref{s-genrank2}(i) and the finiteness of
$\roots (\chi ;\Ndbasis _3,\Ndbasis _1+\Ndbasis _2)$
\Dchaintwo{}{$qr^{-1}$}{$q^{-1}$}{$qrs$}~
imply that
$qrs=-1$ or $rs=1$ or $s=-1,qrs\in R_3$, or $s=-1,qr^{-1}\in R_3$.

If $rs\not=1$ then $s=-1$ and it remains to consider the case $-qr\in R_3$.
However Proposition \ref{s-genrank2}(i) and the finiteness of
$\roots (\chi ;\Ndbasis _1+\Ndbasis _2+\Ndbasis _3,\Ndbasis _1)$
\Dchaintwo{}{$-q$}{$qr$}{$q$}~ give that
$qr=1$ or $q^4r^2=1$. The latter equation is a contradiction
to $(-qr)^3=1$ and $qr^{-1}\not=-1$.

If $rs=1$ then the finiteness of
$\roots (\chi ;\Ndbasis _1+\Ndbasis _2+\Ndbasis _3,\Ndbasis _1)$
\Dchaintwo{}{$qr^{-1}$}{$qr$}{$q$}~ implies that
$(qr-1)(q+1)(q^2r-1)(q^2-1)=0$ or $q^4r=-1$, $\{q,qr^{-1}\}\cap R_3\not=\{\}$.
If $q^2r=1$ then relation $(qr^{-1})^2\not=1$ gives that $q^6\not=1$.
By Theorem \ref{t-Cartan} $(\chi ,E)$ is not of Cartan type. Thus
the finiteness of $\roots (\chi ;\Ndbasis _1,\Ndbasis _3)$ yields
$q^3\in R_3\cup R_4\cup R_5$. Using this
the finiteness of
$\roots (\chi ;\Ndbasis _1+2\Ndbasis _3,\Ndbasis _1+\Ndbasis _2)$
\Dchaintwo{}{$q^{11}$}{$q^{-2}$}{$q$}~ implies that $qr^{-1}=q^3\in R_3$.

Finally, if relations $rs=1$, $q^4r=-1$, and $q\in R_3$ hold,
then the finiteness of
$\roots (\chi ;\Ndbasis _1+\Ndbasis _2+\Ndbasis _3,\Ndbasis _1)$
\Dchaintwo{}{$-q^{-1}$}{$-1$}{$q$}~ gives a contradiction to Lemma
\ref{l-specrank2}(i).
\end{bew}

\begin{lemma}\label{l-a2leg}
If $(\roots ,\chi ,E)$ is a connected arithmetic root system with generalized
Dynkin diagram
\begin{align*}
\Dchainthree{}{$q$}{$q^{-1}$}{$q$}{$r$}{$s$}
\end{align*}
such that $r,s\not=1$, $q^2\not=1$, then $(qr-1)(q^2r-1)(q^3r-1)(r^4s^2-1)=0$.
\end{lemma}

\begin{bew}
Assume that the vertices of $\gDd _{\chi ,E}$ correspond from left to right to
$\Ndbasis _3,\Ndbasis _1$ and $\Ndbasis _2$, respectively.
If $qr\not=1$ then $2\Ndbasis _1+\Ndbasis _2\in \roots $ by Lemma
\ref{l-posroots}(i). Further, $2\Ndbasis _1+\Ndbasis _3\notin \roots $.
Since $4\Ndbasis _1+2\Ndbasis _2+2\Ndbasis _3\notin \roots $ by Lemma
\ref{l-multroots}, one has either $4\Ndbasis _1+\Ndbasis _2+2\Ndbasis _3
\notin \roots $ or $\chi \chi \op (4\Ndbasis _1+\Ndbasis _2+2\Ndbasis _3,
\Ndbasis _2)=1$. In the first case Lemma \ref{l-posroots}(i) applied to
$\roots (\chi ;\Ndbasis _1+\Ndbasis _3,2\Ndbasis _1+\Ndbasis _2)$
\Dchaintwo{}{$q$}{$q^2r$}{$q^4r^2s$}~ implies that $q^2r=1$ or $q^3r=1$,
and in the second one gets $r^4s^2=1$.
\end{bew}

\begin{lemma}\label{l-Weyleqtopath}
If $(\roots ,\chi ,E)$ is a connected arithmetic root system such that
$(q_{12}q_{21}-1)(q_{13}q_{31}-1)(q_{23}q_{32}-1)\not=0$ then
it is Weyl equivalent to an arithmetic root system
$(\roots ',\chi ',E)$ such that $q'_{11}=-1$ and $q'_{23}q'_{32}=1$,
where $q'_{ij}=\chi '(\Ndbasis _i,\Ndbasis _j)$.
\end{lemma}

\begin{bew}
By Lemma \ref{l-cycle} one can assume that $q_{11}=-1$.
Set $\chi '(\Ndbasis _i,\Ndbasis _j):=
\chi (s_{\Ndbasis _1,E}(\Ndbasis _i),s_{\Ndbasis _1,E}(\Ndbasis _j))$
for all $i,j\in \{1,2,3\}$ and apply
Example \ref{ex-Weyleq} with $n=3$
and $i=1$. One obtains $p_{12}=-q_{12}q_{21}$, $p_{13}=-q_{13}q_{31}$,
$q'_{11}=-1$,
and $q'_{23}q'_{32}=p_{12}p_{13}q_{23}q_{32}=1$ by Lemma
\ref{l-trianglecond}.
\end{bew}

\begin{lemma}\label{l-Weyleqtopath2}
If $q_{12}q_{21}\not=1$, $q_{13}q_{31}\not=1$, and $q_{23}q_{32}=1$
then $(\roots ,\chi ,E)$ is Weyl equivalent to an arithmetic root system
$(\roots ',\chi ',E)$ such that equations
$(q'_{11}+1)(q'_{11}q'_{13}q'_{31}-1)=0$, $q'_{23}q'_{32}=1$ hold, where
$q'_{ij}=\chi '(\Ndbasis _i,\Ndbasis _j)$.
\end{lemma}

\begin{bew}
By Lemma \ref{l-2forkcond} and twist equivalence it is sufficient to prove the claim
under the condition that all of the relations
$q_{12}q_{21}\not=1$, $q_{13}q_{31}\not=1$,
$q_{23}q_{32}=1$, $q_{11}q_{12}q_{21}q_{13}q_{31}=-1$,
$q_{11}^2\not=1$, $q_{11}q_{12}q_{21}\not=1$, and
$q_{11}q_{13}q_{31}\not=1$ hold.

Step 1. \textit{If $q_{22}=-1$ or $q_{33}=-1$ then the claim holds.}
By twist equivalence one can assume that $q_{22}=-1$.
Then Example \ref{ex-Weyleq} with $n=3$ and $i=2$ gives that
$q'_{22}=-1$, $q'_{12}q'_{21}=q_{12}^{-1}q_{21}^{-1}$, $q'_{11}=-q_{11}q_{12}q_{21}$,
$q'_{13}q'_{31}=q_{13}q_{31}$, $q'_{33}=q_{33}$, $q'_{23}q'_{32}=1$, and hence
$q'_{11}q'_{13}q'_{31}=1$.

Step 2. \textit{If $q_{13}q_{31}q_{33}\not=1$ then the lemma is true.}
By Step 1 one can assume that $q_{33}\not=-1$. Since $q_{11}\not=-1$
and $q_{11}q_{13}q_{31}\not=1$, Proposition \ref{s-genrank2}(i)
gives that $q_{11}q_{13}^2q_{31}^2q_{33}=-1$ and one of the relations
$q_{11}\in R_3$, $q_{33}\in R_3$ is valid.

Suppose first that $q_{11}\in R_3$. Then
Example \ref{ex-Weyleq} with $n=3$ and $i=1$ together with equation
$q_{11}q_{12}q_{21}q_{13}q_{31}=-1$ gives that
$p_{13}=q_{11}^{-1}q_{13}q_{31}$, $q'_{33}=p_{13}^2q_{33}=-1$, and
$q'_{23}q'_{32}=1$, and hence the lemma
follows from Step 1. On the other hand, if $q_{33}\in R_3$ and
$q_{11}q_{13}^2q_{31}^2q_{33}=-1$
then Example \ref{ex-Weyleq} with $n=3$ and $i=3$ gives
$p_{31}=q_{33}^{-1}q_{13}q_{31}$ and $q'_{11}=p_{31}^2q_{11}=-1$.

Step 3. \textit{If $q_{12}q_{21}q_{22}=1$, $q_{13}q_{31}q_{33}=1$,
and $(q_{11}^2+q_{11}+1)(q_{11}^2q_{12}q_{21}-1)=0$ then the assertion of
the lemma holds.}
Assume first that $q_{11}\in R_3$. Apply Example \ref{ex-Weyleq} with $i=1$
and note that $q_{11}q_{12}q_{21}q_{13}q_{31}=-1$, $q_{11}q_{12}q_{21}\not=1$,
$q_{11}q_{13}q_{31}\not=1$. Thus $q'_{11}=q_{11}$,
$q'_{12}q'_{21}={q'}_{22}^{-1}=q_{11}^{-1}q_{12}^{-1}q_{21}^{-1}$,
$q'_{13}q'_{31}={q'}_{33}^{-1}=q_{11}^{-1}q_{13}^{-1}q_{31}^{-1}$,
and $q'_{23}q'_{32}=1$.
Thus $W_{\chi ,E}=W_{\chi '',E}$, where the structure constants
$q''_{ij}$ of $(\chi '',E)$ satisfy the relations
${q''}_{11}^2q''_{12}q''_{21}={q''}_{11}^2q''_{13}q''_{31}=
q''_{12}q''_{21}q''_{22}=q''_{13}q''_{31}q''_{33}=q''_{23}q''_{32}=1$
and ${q''}_{11}^2\not=1$,
that is $(\chi '',E)$ is of infinite Cartan type. By Theorem \ref{t-Cartan}
this is a contradiction
to the finiteness of $\roots $.

Suppose now that $q_{11}\notin R_3$, and hence $q_{11}^2q_{12}q_{21}=1$.
If $q_{11}^aq_{13}q_{31}=1$ for some $a\ge 2$ then $(\chi ,E)$ is of infinite
Cartan type, which is a contradiction. On the other hand, if
$q_{11}\in R_{m+1}$ with $m\ge 3$, and $(q_{13}q_{31})^{m+1}\not=1$, then
in Example \ref{ex-Weyleq} with $i=1$ one gets
$q'_{23}q'_{32}=(q_{11}^{-1}q_{13}q_{31})^2$.
If $q'_{23}q'_{32}\not=1$ then the Lemma holds by Lemma \ref{l-Weyleqtopath}.
If $q'_{23}q'_{32}=1$ then $q_{13}q_{31}=-q_{11}$ which contradicts the
finiteness of $\roots (\chi ;\Ndbasis _1,\Ndbasis _3)$, Lemma
\ref{l-specrank2}(iii), and relation $(q_{13}q_{31})^{m+1}\not=1$.

Now turn to the proof of the lemma.
By Lemma \ref{l-t2t3t7t14} one can assume that
$(2a+1)\Ndbasis _1+2\Ndbasis _2\notin \roots $ for all $a\in \ndN _0$.
If $q_{22}=-1$ then the claim of the lemma holds by Step 1.
Otherwise $q_{12}q_{21}q_{22}=1$ and either $q_{11}\in R_3$ or
$q_{11}^2q_{12}q_{21}=1$. Hence Steps 2 and 3 prove the lemma.
\end{bew}

\begin{lemma}\label{l-Weyleqtopath3}
If $(\roots ,\chi ,E)$ is a connected arithmetic root system 
then either $\gDd _{\chi ,E}$ appears in rows 9--11 of Table 2, or
it is Weyl equivalent to an arithmetic root system
$(\roots ',\chi ',E)$ such that $q'_{23}q'_{32}=1$ and $q'_{11}q'_{13}q'_{31}=1$,
where $q'_{ij}=\chi '(\Ndbasis _i,\Ndbasis _j)$.
\end{lemma}

\begin{bew}
By Lemmata \ref{l-Weyleqtopath} and \ref{l-Weyleqtopath2} it is sufficient to prove
the claim under the additional assumption that $q_{11}=-1$ and $q_{23}q_{32}=1$,
where $q_{ij}=\chi (\Ndbasis _i,\Ndbasis _j)$.
By Lemma \ref{l-t2t3t7t14} one can assume that $q_{22}=-1$
or $q_{12}q_{21}q_{22}=1$. In the first case one can apply Example \ref{ex-Weyleq}
with $n=3$ and $i=2$. One gets $p_{21}=-q_{12}q_{21}$, $m_{21}=1$, $p_{23}=1$,
$m_{23}=0$, and hence $q'_{11}q'_{12}q'_{21}=-q_{12}q_{21}q_{11}
q_{12}^{-2}q_{21}^{-2}q_{12}q_{21}=1$ and $q'_{23}q'_{32}=1$, and therefore the
assertion follows using twist equivalence.

Assume now that $q_{11}=-1$, $q_{22}\not=-1$, $q_{33}\not=-1$,
and $q_{12}q_{21}q_{22}=1$. If
$q_{13}q_{31}q_{33}=1$ then it can happen that $q_{12}q_{21}q_{13}q_{31}=1$,
and hence one has an example in row 8 of Table 2, and the claim of
the lemma holds. Otherwise $\gDd _{\chi ,E}$ appears in rows 9--11 of Table 2.
Therefore one can assume additionally that $q_{13}q_{31}q_{33}\not=1$.
Then one can apply Lemma \ref{l-bigleg}
to conclude that either $q_{12}q_{21}q_{13}q_{31}=1$ or
$q_{12}q_{21}q_{13}^2q_{31}^2q_{33}=-1$.
If $q_{12}q_{21}q_{13}q_{31}=1$ then using Example \ref{ex-Weyleq} with
$n=3$ and $i=1$ one gets $q'_{23}q'_{32}=1$, $q'_{11}=q'_{22}=-1$, and hence
the claim of the lemma holds by the first part of the proof. On the other hand,
if $q_{12}q_{21}q_{13}^2q_{31}^2q_{33}=-1$ and $q_{12}q_{21}q_{13}q_{31}\not=1$
then one can use twice Example \ref{ex-Weyleq}, first with $i=1$ and then with
$i=2$, to obtain a Weyl equivalent arithmetic root system
$(\roots ',\chi ',E)$ such that $q'_{22}=q'_{33}=-1$, $q'_{23}q'_{32}\not=1$,
and $q'_{13}q'_{31}=1$, and hence the assertion
of the lemma holds again by the first paragraph of the proof.
\end{bew}

%

In what follows Theorem \ref{t-class3} will be reformulated in two pieces
and the proofs will be given.

\begin{thm}\label{t-finiteness}
If $\gDd _{\chi ,E}$ appears in Table 2
then $(\roots ,\chi ,E)$ is an arithmetic root system, where
$\roots =\bigcup \{F\,|\,(\id ,F)\in W_{\chi ,E}\}$.
Further, two such arithmetic root systems
$(\roots ,\chi ,E)$, $(\roots ',\chi ',E)$
are Weyl equivalent if and only if $\gDd _{\chi ,E}$ and $\gDd _{\chi ',E}$
appear in the same row of Table 2 and can be presented with the same set of
fixed parameters.
\end{thm}

\begin{bew}
The second assertion can be checked case by case using Example \ref{ex-Weyleq}.
Therefore $\grgDd ^{\chi ,E}$ is finite for all pairs $(\chi, E)$ such that
$\gDd _{\chi ,E}$ appears in Table 2. Now the first assertion follows
by direct calculations from
Lemmata \ref{l-genWBg} and \ref{l-finextWBG}. The structure of
$\WBg _{\chi ,E}$ is given in Table 2. Here $[m]$ denotes the Coxeter group
generated by two elements $s_1,s_2$ and the relations
$s_1^2=s_2^2=(s_1s_2)^m=1$,
and $[m_1,m_2]$ is the Coxeter group generated by three elements
$s_1,s_2,s_3$ and the relations
$s_1^2=s_2^2=s_3^2=(s_1s_2)^{m_1}=(s_1s_3)^2=(s_2s_3)^{m_2}=1$.
\end{bew}

\begin{thm}\label{t-getTable2}
If $(\roots ,\chi ,E)$ is a connected arithmetic root system then
$\gDd _{\chi ,E}$ appears in Table 2.
\end{thm}

\begin{bew}
By Theorem \ref{t-finiteness} it suffices to consider $(\roots ,\chi ,E)$
up to Weyl equivalence. By Lemma \ref{l-Weyleqtopath3} one can assume that
$q_{23}q_{32}=1$ and $q_{11}q_{13}q_{31}=1$. Further,
by Lemma \ref{l-bigleg} one can suppose that equation
$(q_{33}+1)(q_{33}-q_{11})(q_{11}-q_{12}q_{21}q_{33})
(q_{11}q_{12}q_{21}-1)=0$ holds.

Step 1. \textit{Up to Weyl equivalence one can assume that
$q_{11}q_{13}q_{31}=1$ and $(q_{33}+1)(q_{33}-q_{11})=0$.}
Assume that $(q_{33}+1)(q_{11}-q_{33})\not=0$.
If $q_{11}q_{12}q_{21}=1$ then by Lemma \ref{l-t2t3t7t14}
one has $(q_{22}+1)(q_{12}q_{21}q_{22}-1)=0$
and hence the claim of Step 1 holds by exchanging
$\Ndbasis _2$ and $\Ndbasis _3$.
Otherwise the last equation above Step 1
implies that $q_{11}=q_{12}q_{21}q_{33}$. Further, by Lemma
\ref{l-specbigleg} one obtains that either $q_{11}=-1$ or
$q_{33}\in R_3$. In the second case ($q_{33}\in R_3$)
one can apply Example \ref{ex-Weyleq} with $i=3$ ($q:=q_{11}$,
$r:=q_{12}q_{21}$, $s:=q_{22}$). Here and in the following
the node corresponding to the basis vector $\Ndbasis _i$
and the reflection $s_{\Ndbasis _i,E}$, respectively,
will be denoted by a black circle.
\begin{gather*}
\Dchainthree{l}{$qr^{-1}$}{$q^{-1}$}{$q$}{$r$}{$s$} ~ \Rightarrow ~
\Dchainthree{}{$qr^{-1}$}{$r$}{$r^{-1}$}{$r$}{$s$}.
\end{gather*}
Using additionally
equation $(q_{22}+1)(q_{12}q_{21}q_{22}-1)=0$, which follows from Lemma
\ref{l-t2t3t7t14}, one obtains an
arithmetic root system which satisfies the conditions of Step 1 up to twist
equivalence. In the final case $q_{11}=-1$ one has again the equation
$(q_{22}+1)(q_{12}q_{21}q_{22}-1)=0$ from Lemma \ref{l-t2t3t7t14}.
Thus the claim follows from the following transformations
($r=q_{12}q_{21}$).
\begin{gather*}
\Dchainthree{r}{$-r^{-1}$}{$-1$}{$-1$}{$r$}{$-1$} ~\Rightarrow ~
\Dchainthree{}{$-r^{-1}$}{$-1$}{$r$}{$r^{-1}$}{$-1$}~,\\
\rule{0pt}{10\unitlength}
\Dchainthree{m}{$-r^{-1}$}{$-1$}{$-1$}{$r$}{$r^{-1}$} ~\Rightarrow ~
\cDtriangle{r}{$-r^{-1}$}{$-1$}{$-1$}{$-1$}{$r^{-1}$}{$-r$} ~\Rightarrow ~
\Dchainthree{r}{$r^{-1}$}{$r$}{$-1$}{$-r^{-1}$}{$-1$} ~\Rightarrow ~
\Dchainthree{}{$r^{-1}$}{$r$}{$-r^{-1}$}{$-r$}{$-1$}~.
\end{gather*}

Using Step 1 and Proposition \ref{s-genrank2}(i) for $\roots (\chi ;
\Ndbasis _1,\Ndbasis _2)$, 9 special cases will be considered which
together prove the claim of the theorem. In all of them the vertices of
$\gDd _{\chi ,E}$ correspond from left to right to the basis vectors
$\Ndbasis _3,\Ndbasis _1$, and $\Ndbasis _2$ of $E$, respectively.

Step 2. \textit{If $q_{11}=q_{13}q_{31}=q_{33}=-1$ then the claim of the
theorem holds.} Exchanging $\Ndbasis _2$ and $\Ndbasis _3$,
Lemma \ref{l-bigleg} implies that $q_{22}=-1$ or $q_{12}q_{21}q_{22}=1$ or
$q_{12}^2q_{21}^2q_{22}=1$ or $q_{12}q_{21}=-1$. If $q_{12}q_{21}q_{22}=1$
then $\gDd _{\chi ,E}$ appears in one of rows 1,9, or 10 of Table 2.
If $q_{22}=-1$ then the finiteness of
$\roots (\chi ;\Ndbasis _3,\Ndbasis _1+\Ndbasis _2)$
\Dchaintwo{}{$-1$}{$-1$}{~$q_{12}q_{21}$}~~ and Lemma \ref{l-specrank2}(i)
imply that $q_{12}q_{21}\in R_2\cup R_3\cup R_4\cup R_6$.
Then $\gDd _{\chi ,E}$ appears in row 1,17,6 and 7 of Table 2, respectively.
If $q_{12}^2q_{21}^2q_{22}=1$ then by the transformations ($r=q_{12}q_{21}$)
\begin{gather*}
\Dchainthree{m}{$-1$}{$-1$}{$-1$}{$r$}{$r^{-2}$} ~\Rightarrow ~
\rule{0pt}{10\unitlength}
\cDtriangle{l}{$-1$}{$-1$}{$-r^{-1}$}{$-1$}{$r^{-1}$}{$-r$} ~\Rightarrow ~
\Dchainthree{}{$-1$}{$-r^{-1}$}{$-1$}{$-1$}{$-1$}
\end{gather*}
all such arithmetic root systems are Weyl equivalent to one of the previous
case.

Finally, if $q_{12}q_{21}=-1$ then the finiteness of
$\roots (\chi ;\Ndbasis _1,\Ndbasis _2)$ and Lemma \ref{l-specrank2}(i) give
that $q_{22}\in R_2\cup R_3\cup R_4\cup R_6$. If $q_{22}\in R_2\cup R_4$
then $\gDd _{\chi ,E}$ appears in rows 1 and 2 of Table 2, respectively.
If $q_{22}\in R_6$ then $(\chi ,E)$ is of infinite Cartan type which contradicts
Theorem \ref{t-Cartan}. If $q_{22}\in R_3$ then
$\Ndbasis _1+2\Ndbasis _2+\Ndbasis _3\in
\roots (\chi ;\Ndbasis _1+\Ndbasis _3,\Ndbasis _2)$
\Dchaintwo{}{$-1$}{$-1$}{$q_{22}$}
and hence the finiteness of
$\roots (\chi ;\Ndbasis _1+2\Ndbasis _2+\Ndbasis _3,\Ndbasis _1+\Ndbasis _2)$
\Dchaintwo{}{$-q_{22}$}{$q_{22}$}{$q_{22}$}
contradicts Proposition \ref{s-genrank2}(i).

Step 3. \textit{If $\gDd _{\chi ,E}$ is
\Dchainthree{}{$-1$}{$q^{-1}$}{$q$}{$r$}{$-1$}~ with $q^2\not=1$, $r\not=1$,
then the claim of the theorem holds.} The finiteness of
$\roots (\chi ;\Ndbasis _1+\Ndbasis _2+\Ndbasis _3,\Ndbasis _1)$
\Dchaintwo{}{$r$}{$qr$}{$q$} and Proposition \ref{s-genrank2}(i) imply that
$(qr-1)(q^2r-1)(r+1)(qr^2-1)=0$ or $q^3r^3=-1$, $\{q,r\}\cap R_3\not=\{\}$.
However if $(qr-1)(q^2r-1)(r+1)(qr^2-1)\not=0$ and say $q\in R_3$ then
by Table 1 one has $r\in R_3\cup R_4\cup R_{18}\cup R_{30}$ which is a
contradiction to $r^3=-1$.

In the case $qr=1$ $\gDd _{\chi ,E}$ appears in row 8 of Table 2.

If $q^2r=1$ and $r\not=-1$ then the finiteness of
$\roots (\chi ;\Ndbasis _1+\Ndbasis _2,\Ndbasis _3)$
\Dchaintwo{}{$-q^{-1}$}{$q^{-1}$}{$-1$}~ and Lemma \ref{l-specrank2}(ii)
imply that $q\in R_3\cup R_5\cup R_6\cup R_8$. If $q\in R_3\cup R_6$
then $\gDd _{\chi ,E}$ appears in row 15 and 17, respectively.
On the other hand, if $q\in R_5$
then $\Ndbasis _1+2\Ndbasis _2+\Ndbasis _3\notin
\roots (\chi ;\Ndbasis _1+\Ndbasis _3,\Ndbasis _2)$
\Dchaintwo{}{$-1$}{$q^{-2}$}{$-1$}~,
but $2\Ndbasis _1+2\Ndbasis _2+\Ndbasis _3\in
\roots (\chi ;\Ndbasis _1+\Ndbasis _2,\Ndbasis _3)$ by Lemma
\ref{l-posroots}(i). However $\roots (\chi ;
2\Ndbasis _1+2\Ndbasis _2+\Ndbasis _3,\Ndbasis _1)$
\Dchaintwo{}{$-q^{-1}$}{$q^{-1}$}{$q$}~ is of infinite Cartan type
which is a contradiction to Theorem \ref{t-Cartan} and Proposition
\ref{s-subARS}. Finally, if $q\in R_8$ then
$4\Ndbasis _1+3\Ndbasis _2+3\Ndbasis _3\in
\roots (\chi;\Ndbasis _1+\Ndbasis _2+\Ndbasis _3,\Ndbasis _1)$
\Dchaintwo{}{$q^{-2}$}{$q^{-1}$}{$q$}~ and
$4\Ndbasis _1+3\Ndbasis _2+\Ndbasis _3\in
\roots (\chi;\Ndbasis _1+\Ndbasis _2,\Ndbasis _1+\Ndbasis _3)$
\Dchaintwo{}{$-q^{-1}$}{$q^{-1}$}{$-1$}~ which is a contradiction
to $\roots (\chi ;4\Ndbasis _1+3\Ndbasis _2+\Ndbasis _3,\Ndbasis _3)=
\roots \cap \linc{\ndZ }{\{\Ndbasis _3,
4\Ndbasis _1+3\Ndbasis _2+\Ndbasis _3\}}$
\Dchaintwo{}{$-1$}{$-1$}{$-1$}.

It remains to consider the cases $r=-1$ and $qr^2=1$, respectively.
Here one can use the transformations
\begin{gather*}
\Dchainthree{l}{$-1$}{$q^{-1}$}{$q$}{$-1$}{$-1$}~\Rightarrow~
\Dchainthree{}{$-1$}{$q$}{$-1$}{$-1$}{$-1$}~,\\
\Dchainthree{l}{$-1$}{$r^2$}{$r^{-2}$}{$r$}{$-1$}~\Rightarrow~
\Dchainthree{r}{$-1$}{$r^{-2}$}{$-1$}{$r$}{$-1$}~\Rightarrow~
\Dchainthree{}{$-1$}{$r^{-2}$}{$r$}{$r^{-1}$}{$-1$}~
\end{gather*}
to obtain arithmetic root systems which were already considered in Step 2
and Step 3, respectively.

Step 4. \textit{If $\gDd _{\chi ,E}$ is
\Dchainthree{}{$-1$}{$q^{-1}$}{$q$}{$q^{-1}$}{$s$}~ with $q^2\not=1$,
$s^2\not=1$,
then the claim of the theorem holds.}
Since $\roots (\chi ;\Ndbasis _1,\Ndbasis _2)$
\Dchaintwo{}{$q$}{$q^{-1}$}{$s$}
and $\roots (\chi ;\Ndbasis _1+\Ndbasis _3,\Ndbasis _2)$
\Dchaintwo{}{$-1$}{$q^{-1}$}{$s$}
are finite, Lemma \ref{l-specrank2-2} implies that either $s^m=q$
for $m\in \{1,2,3\}$ or $s\in R_3$. If $s=q$ or $s^2=q$ then
$\gDd _{\chi ,E}$ appears in row 4 and row 5 of Table 2, respectively.
If $s^3=q$ then Proposition
\ref{s-genrank2}(i) and the finiteness of
$\roots (\chi ;\Ndbasis _1+\Ndbasis _2+\Ndbasis _3,\Ndbasis _1+2\Ndbasis _2)$
\Dchaintwo{}{$-s^{-2}$}{$s^{-2}$}{$s$}~
give that either $s^4=-1$ or $s^5=1$, $\{s,-s^{-2}\}\cap R_3\not=\{\}$.
The second case contradicts himself. If $s\in R_8$ then
one has $2\Ndbasis _1+3\Ndbasis _2\in \roots $ and
$2\Ndbasis _1+3\Ndbasis _2+2\Ndbasis _3\notin \roots (\chi ;
2\Ndbasis _1+3\Ndbasis _2,\Ndbasis _3)$
\Dchaintwo{}{$s^3$}{$s^{-6}$}{$-1$}~.This contradicts to the fact
that $2(\Ndbasis _1+\Ndbasis _3)+3\Ndbasis _2\in
\roots (\chi ;\Ndbasis _1+\Ndbasis _3,\Ndbasis _2)$.

Suppose now that $s\in R_3$ and $q^3\not=1$.
Then Proposition \ref{s-genrank2}(i) and the finiteness of
$\roots (\chi ;\Ndbasis _1+\Ndbasis _2+\Ndbasis _3,\Ndbasis _1+2\Ndbasis _2)$
\longDchaintwo{}{$-q^{-1}s$}{$q^{-2}s$}{$q^{-1}s$}~
imply that equation $(q^2-s)(q^{-3}s^2-1)(q^{-3}s^2+1)=0$ holds.
If $s=q^2$ then $q^3\not=1$ gives that $q\in R_6$, $s=-q^{-1}$. In this case
$\gDd _{\chi ,E}$ appears in row 14 of Table 2.
On the other hand, if $s^2=q^3$ then $s=q^{-3}$, $q\in R_9$, and
$\roots (\chi ;\Ndbasis _1+\Ndbasis _2+\Ndbasis _3,\Ndbasis _1+2\Ndbasis _2)$
corresponds to the generalized Dynkin diagram
\Dchaintwo{}{$-q^{-4}$}{$q^{-5}$}{$q^{-4}$}~. This contradicts
Theorem \ref{t-Cartan}. In the same way $s^2=-q^3$ leads
to a contradiction.

Step 5. \textit{If $\gDd _{\chi ,E}$ is
\Dchainthree{}{$-1$}{$q^{-1}$}{$q$}{$r$}{$r^{-1}$}~ with $q^2\not=1$,
$r^2\not=1$, $qr\not=1$, then the claim of the theorem holds.}
By Table 1 and the finiteness of $\roots (\chi ;\Ndbasis _1,\Ndbasis _2)$
\Dchaintwo{}{$q$}{$r$}{$r^{-1}$}~
one has either $q^mr=1$ with $m\in \{2,3\}$ or $q\in R_{m+1}$ with
$m\in \{2,3,4\}$ and $r^{m+1}\not=1$.
If $q^2r=1$ or $q^3r=1$ then $\gDd _{\chi ,E}$ appears in row
6 and 7 of Table 2, respectively. If $q\in R_5$ then $r\in R_{30}$
and $q=-r^{-3}$ again
by the finiteness of $\roots (\chi ;\Ndbasis _1,\Ndbasis _2)$.
Then $\roots (\chi ;2\Ndbasis _1+\Ndbasis _2,
\Ndbasis _3)$
\Dchaintwo{}{$-r^4$}{$r^6$}{$-1$}~ is not finite.

If $q\in R_4$ then by the finiteness of $\roots (\chi ;\Ndbasis _1,\Ndbasis _2)$
one gets either $r\in R_{24}$, $q=r^6$, or $r\in R_8$, $q=r^2$.
Since $\chi \chi \op (2\Ndbasis _1+\Ndbasis _2,\Ndbasis _1+\Ndbasis _3)=
-r\not=1$, one has $3\Ndbasis _1+\Ndbasis _2+\Ndbasis _3\in \roots $.
In both cases $\roots (\chi ;3\Ndbasis _1+\Ndbasis _2+\Ndbasis _3,
\Ndbasis _1+\Ndbasis _2)$
\Dchaintwo{}{$r^2$}{$qr^2$}{$q$}~ is not finite.

Finally, if $q\in R_3$ then using Table 1 (or Theorem 5 in
\cite{a-Heck04b}) the finiteness of
$\roots (\chi ;2\Ndbasis _1+\Ndbasis _2,\Ndbasis _3)$
\Dchaintwo{}{$qr$}{$q$}{$-1$}~ implies that either $q^2r=1$
or $qr=-1$ or $(qr)^2q=1$. The case $q^2r=1$ was already considered,
and $r^2=1$ was excluded.
Thus $qr=-1$ and then $\gDd _{\chi ,E}$ appears in row 16 of Table 2.

Step 6. \textit{If $\gDd _{\chi ,E}$ is
\Dchainthree{}{$-1$}{$q^{-1}$}{$q$}{$r$}{$s$}~ with $q^2\not=1$,
$s^2\not=1$, $qr\not=1$, and $rs\not=1$, then the claim of the theorem holds.}
By Proposition \ref{s-genrank2}(i) one has $qr^2s=-1$ and $\{q,s\}\cap R_3
\not=\{\}$.
If $q\in R_3$ then $r\not=q$. The transformations
\begin{align*}
\rule{0pt}{10\unitlength}
\Dchainthree{m}{$-1$}{$q^{-1}$}{$q$}{$r$}{$s$} ~\Rightarrow ~
\cDtriangle{l}{$-1$}{$q$}{$-1$}{$q^{-1}$}{$q^{-1}r^{-1}$}{$q^{-1}r$}
~~~\Rightarrow ~
\Dchainthree{l}{$-1$}{$q$}{$-1$}{$qr^{-1}$}{~~$q^{-1}r$} ~\Rightarrow ~
\Dchainthree{}{$-1$}{$q^{-1}$}{$q$}{$qr^{-1}$}{~~$q^{-1}r$}~
\end{align*}
give that $(\roots ,\chi ,E)$ is Weyl equivalent to an arithmetic root system
which was already considered in the previous steps.

When $s\in R_3$ and $q\notin R_3$ then by Table 1 there are only three
possibilities: either $r\in R_{24}$, $q=r^{-6}$, or
$q\in R_{18}$, $r=q^{-2}$, or $q\in R_{30}$, $r=q^{-3}$.
In all cases
$\roots (\chi ;\Ndbasis _1+\Ndbasis _2,\Ndbasis _1+\Ndbasis _3)$
\Dchaintwo{}{$-r^{-1}$}{$qr$}{$-1$}~ is not finite by Table 1.

Step 7. \textit{If $\gDd _{\chi ,E}$ is
\Dchainthree{}{$q$}{$q^{-1}$}{$q$}{$r$}{$-1$}~ with $q^2\not=1$, $r\not=1$,
then the claim of the theorem holds.}
By Lemma \ref{l-a2leg} one has $(qr-1)(q^2r-1)(q^3r-1)(r^4-1)=0$.
If $qr=1$ then $\gDd _{\chi ,E}$ appears in row 4 of Table 2.
If $q^2r=1$ then the finiteness of $\roots (\chi ;\Ndbasis _3,
\Ndbasis _1+\Ndbasis _2)$
\Dchaintwo{}{$q$}{$q^{-1}$}{$-q^{-1}$}~
and Lemma \ref{l-specrank2}(iii) imply that $q\in R_3\cup R_4\cup R_6\cup R_8$.
If $q\in R_3\cup R_6\cup R_4$ then $\gDd _{\chi ,E}$ appears in row
13 and 3 of Table 2, respectively.
If $q\in R_8$ then $2\Ndbasis _1+3\Ndbasis _2+2\Ndbasis _3\notin
\roots (\chi ;\Ndbasis _1+\Ndbasis _3,\Ndbasis _2)$
\Dchaintwo{}{$q$}{$q^{-2}$}{$-1$}~ and
$3\Ndbasis _1+3\Ndbasis _2+2\Ndbasis _3\in
\roots (\chi ;\Ndbasis _1+\Ndbasis _2,\Ndbasis _3)$
\Dchaintwo{}{$-q^{-1}$}{$q^{-1}$}{$q$}~, but
$\roots (\chi ;3\Ndbasis _1+3\Ndbasis _2+2\Ndbasis _3,\Ndbasis _1)$
\Dchaintwo{}{$q$}{$q^{-2}$}{$q$}~ is not finite.

Suppose now that $q^3r=1$. Then Proposition \ref{s-genrank2}(i) and
the finiteness of $\roots (\chi ;\Ndbasis _1+\Ndbasis _2,
\Ndbasis _1+\Ndbasis _3)$ \Dchaintwo{}{$-q^{-2}$}{$q^{-2}$}{$q$} give that
$q\in R_8$. In this case $\roots (\chi ;\Ndbasis _1+\Ndbasis _2,\Ndbasis _3)$
\Dchaintwo{}{$-q^{-2}$}{$q^{-1}$}{$q$}~ is not finite.

Finally, if $r^4=1$ and $q^mr\not=1$ for all $m\in \{1,2,3\}$
then $r\in R_2\cup R_4$. The finiteness of
$\roots (\chi ;\Ndbasis _3,2\Ndbasis _1+\Ndbasis _2)$
\Dchaintwo{}{$q$}{$q^{-2}$}{$-q^4r^2$} and Proposition \ref{s-genrank2}(i)
imply that either $r=-1$, $q^4=1$, or $r\in R_4$, $q^4=-1$.
Since $q^2r\not=1$, only the case $q\in R_8$, $q^6r=1$ can appear,
but this is a contradiction to the finiteness of
$\roots (\chi ;\Ndbasis _1,\Ndbasis _2)$.

Step 8. \textit{If $\gDd _{\chi ,E}$ is
\Dchainthree{}{$q$}{$q^{-1}$}{$q$}{$r$}{$r^{-1}$}~ with $q^2\not=1$,
$r^2\not=1$,
then the claim of the theorem holds.} By Lemma \ref{l-a2leg}
one has $(qr-1)(q^2r-1)(q^3r-1)=0$.
If $qr=1$ or $q^2r=1$ then $\gDd _{\chi ,E}$ appears in row 1 and row 3
of Table 2, respectively. If $q^3r=1$ then $(\chi ,E)$ is of infinite
Cartan type which is a contradiction to Theorem \ref{t-Cartan}.

Step 9. \textit{If $\gDd _{\chi ,E}$ is
\Dchainthree{}{$q$}{$q^{-1}$}{$q$}{$q^{-1}$}{$s$}~ with $q^2\not=1$,
$s^2\not=1$, and $q\not=s$, then the claim of the theorem holds.}
If $s^mq^{-1}=1$ for some $m\ge 2$ then $(\chi ,E)$ is of Cartan type.
Thus by Theorem \ref{t-Cartan} one has $m=2$ and then
$\gDd _{\chi ,E}$ appears in row 2 of Table 2.

Assume now that $s\in R_{m+1}$ with $m\ge 2$, and $q^{m+1}\not=1$.
Then the finiteness of $\roots (\chi ;\Ndbasis _1,\Ndbasis _2)$
gives that $2\le m\le 4$. Moreover, if $m>2$ then either
$q\in R_{30}$, $s=-q^3$, or
$q\in R_{24}$, $s=-q^6$, or $q\in R_8$, $s=-q^2$. In all cases
one obtains a contradiction to the finiteness of
$\roots (\chi ;\Ndbasis _1+2\Ndbasis _2,\Ndbasis _3)$
\Dchaintwo{}{$s^4q^{-1}$}{$q^{-1}$}{$q$}~.

If $m=2$ then the finiteness of
$\roots (\chi ;\Ndbasis _1+\Ndbasis _2+\Ndbasis _3,\Ndbasis _1+2\Ndbasis _2)$
\Dchaintwo{}{$s$}{$sq^{-2}$}{~$sq^{-1}$}~~ and Proposition \ref{s-genrank2}(i)
imply that
$$(sq^{-2}-1)(s^2q^{-2}-1)(sq^{-1}+1)(s^2q^{-3}-1)(sq^{-5}+1)=0.$$
Since $q^3\not=1$, this yields $(sq+1)(sq^{-1}+1)(sq^3-1)(sq^{-5}+1)=0$.
If $sq=-1$ or $s=-q$ or $sq^3=1$ then $\gDd _{\chi ,E}$ appears in row 12,
16 and 18 of Table 2,
respectively. If $s=-q^5$ then either $q\in R_6$, $sq=-1$, or $q\in R_{30}$.
In the second case $\roots (\chi ;\Ndbasis _1+\Ndbasis _2+\Ndbasis _3,
\Ndbasis _1+2\Ndbasis _2)$ is not finite by Table 1.

Step 10. \textit{If $\gDd _{\chi ,E}$ is
\Dchainthree{}{$q$}{$q^{-1}$}{$q$}{$r$}{$s$}~ with $q^2\not=1$,
$s^2\not=1$, $qr\not=1$, and $rs\not=1$,
then the claim of the theorem holds.} By Proposition
\ref{s-genrank2}(i) one has $qr^2s=-1$ and $\{q,s\}\cap R_3\not=\{\}$.
The finiteness of
$\roots (\chi ;\Ndbasis _1+\Ndbasis _2,\Ndbasis _1+\Ndbasis _3)$
\Dchaintwo{}{$-r^{-1}$}{$qr$}{$q$} and Proposition \ref{s-genrank2}(i)
imply that $(q^2r-1)(q^3r-1)=0$. By the finiteness of
$\roots (\chi ;\Ndbasis _1,\Ndbasis _2)$ if $q^2r=1$ then $q\in R_{18}$
and if $q^3r=1$ then $q\in R_{30}$. Both cases contradict to the
finiteness of
$\roots (\chi ;\Ndbasis _1+\Ndbasis _2,\Ndbasis _1+\Ndbasis _3)$.
\end{bew}

\bibliographystyle{mybib}
\bibliography{quantum}

\end{document}